\documentclass[12pt,a4paper]{article}
\usepackage[utf8]{inputenc}

\usepackage{mathrsfs,amsthm,xcolor,verbatim,bbm,amsmath,amsfonts,
amssymb,nicefrac,enumitem,hyperref,bm,mathtools,xparse,etoolbox,graphicx}
\usepackage[margin=1in]{geometry}
\usepackage[capitalise,sort]{cleveref} 

\usepackage{float}

\crefname{enumi}{item}{items}

\crefname{equation}{}{}
\crefname{subsection}{Subsection}{Subsections}
\crefname{figure}{Figure}{Figures}

\hypersetup{
    colorlinks,
    linkcolor={red!50!black},
    citecolor={green},
    urlcolor={blue!80!black}
}

\theoremstyle{plain}
\newtheorem{theorem}{Theorem} [section]
\newtheorem{lemma}[theorem]{Lemma}
\newtheorem{prop}[theorem]{Proposition}
\newtheorem{cor}[theorem]{Corollary}
\newtheorem{setting}[theorem]{Setting}

\theoremstyle{definition}

\newtheorem{case}{Case}
\newtheorem{remark}[theorem]{Remark}

\crefname{case}{Case}{Cases}

\numberwithin{equation}{section}

\DeclareFontEncoding{LS1}{}{}
\DeclareFontSubstitution{LS1}{stix}{m}{n}
\DeclareMathAlphabet{\mathscr}{LS1}{stixscr}{m}{n}

\renewcommand{\P}{\mathbb{P}}

\newcommand{\R}{\mathbb{R}}
\newcommand{\N}{\mathbb{N}}

\newcommand{\ssum}{\textstyle\sum}
\newcommand{\ssuml}{\textstyle\sum\limits}
\newcommand{\tint}{\textstyle\int}

\newcommand{\with}{\curvearrowleft}

\newcommand{\cA}{\mathcal{A}}
\newcommand{\cB}{\mathcal{B}}

\newcommand{\cF}{\mathcal{F}}
\newcommand{\cG}{\mathcal{G}}

\newcommand{\cL}{\mathcal{L}}
\newcommand{\cM}{\mathcal{M}}
\newcommand{\cN}{\mathcal{N}}

\newcommand{\cT}{\mathcal{T}}

\newcommand{\cV}{\mathcal{V}}

\newcommand{\bfk}{\mathbf{k}}

\newcommand{\fC}{\mathfrak{C}}

\newcommand{\fG}{\mathfrak{G}}

\newcommand{\fL}{\mathfrak{L}}

\newcommand{\fb}{\mathfrak{b}}

\newcommand{\fd}{\mathfrak{d}}

\newcommand{\fw}{\mathfrak{w}}

\renewcommand{\emptyset}{\varnothing}

\DeclarePairedDelimiter{\norm}{\lVert}{\rVert}
\DeclarePairedDelimiter{\abs}{\lvert}{\rvert}
\DeclarePairedDelimiter{\rbr}{(}{)}
\DeclarePairedDelimiter{\br}{[}{]}
\DeclarePairedDelimiter{\cu}{\{}{\}}
\DeclarePairedDelimiter{\spro}{\langle}{\rangle}

\newcommand{\Rect}{\mathfrak{R}}
\renewcommand{\d}{ \mathrm{d}}

\newcommand{\qandq}{\quad\text{and}\quad}
\newcommand{\qqandqq}{\qquad\text{and}\qquad}

\newcommand{\indicator}[1]{\mathbbm{1}_{\smash{#1}}}

\usepackage{tikz}
\usetikzlibrary{matrix,chains,positioning,decorations.pathreplacing,arrows}
\usetikzlibrary{shapes,fadings,arrows.meta}
\tikzset{
	font={\fontsize{9pt}{12}\selectfont}}
\usepackage{adjustbox}

\ExplSyntaxOn

\bool_new:N \g_noteobserve

\NewDocumentCommand{\setnote}{}{
  \bool_gset_true:N \g_noteobserve
}

\NewDocumentCommand{\setobserve}{}{
  \bool_gset_false:N \g_noteobserve
}

\NewDocumentCommand{\nobs}{ o }{
  \IfValueT{#1}{
    \str_if_eq:noTF {note} {#1} {
      \bool_gset_true:N \g_noteobserve
    } {
      \str_if_eq:noTF {Note} {#1} {
        \bool_gset_true:N \g_noteobserve
      } {
        \bool_gset_false:N \g_noteobserve
      }
    }
  }
  \bool_if:nTF { \g_noteobserve } {
    \bool_gset_false:N \g_noteobserve
    note
  } {
    \bool_gset_true:N \g_noteobserve
    observe
  }
  \IfValueF{#1}{~}
}

\NewDocumentCommand{\Nobs}{ o }{
  \IfValueT{#1}{
    \str_if_eq:noTF {note} {#1} {
      \bool_gset_true:N \g_noteobserve
    } {
      \str_if_eq:noTF {Note} {#1} {
        \bool_gset_true:N \g_noteobserve
      } {
        \bool_gset_false:N \g_noteobserve
      }
    }
  }
  \bool_if:nTF { \g_noteobserve } {
    \bool_gset_false:N \g_noteobserve
    Note
  } {
    \bool_gset_true:N \g_noteobserve
    Observe
  }
  \IfValueF{#1}{~}
}

\int_new:N \g_furthermore

\NewDocumentCommand{\Moreover}{ o o }{
  \IfValueT{#1}{
    \str_case:nn {#1} {
      {Furthermore} {\int_set:Nn {\g_furthermore} {0}}
      {Moreover} {\int_set:Nn {\g_furthermore} {1}}
      {In~addition} {\int_set:Nn {\g_furthermore} {2}}
      {note} {\bool_gset_true:N \g_noteobserve}
      {observe} {\bool_gset_false:N \g_noteobserve}
    }
    \IfValueT{#2}{
      \str_case:nn {#2} {
        {Furthermore} {\int_set:Nn {\g_furthermore} {0}}
        {Moreover} {\int_set:Nn {\g_furthermore} {1}}
        {In~addition} {\int_set:Nn {\g_furthermore} {2}}
        {note} {\bool_gset_true:N \g_noteobserve}
        {observe} {\bool_gset_false:N \g_noteobserve}
      }
    }
  }
  \int_case:nn { \int_mod:nn {\g_furthermore} {3} } {
    { 0 } { Furthermore,~\nobs that}
    { 1 } { Moreover,~\nobs that}
    { 2 } { In~addition,~\nobs that}
  }
  \int_incr:N \g_furthermore
  \IfValueF{#1}{~}
}

\bool_new:N \g_hencetherefore

\NewDocumentCommand{\hence}{}{
  \bool_if:nTF { \g_hencetherefore } {
    \bool_gset_false:N \g_hencetherefore
    hence~
  } {
    \bool_gset_true:N \g_hencetherefore
    therefore~
  }
}

\NewDocumentCommand{\Hence}{}{
  \bool_if:nTF { \g_hencetherefore } {
    \bool_gset_false:N \g_hencetherefore
    Hence,~we~obtain~
  } {
    \bool_gset_true:N \g_hencetherefore
    Therefore,~we~obtain~
  }
}

\seq_new:N \g_cflist_loaded
\seq_new:N \g_cflist_pending

\NewDocumentCommand{\cfadd}{ m }
{
	\seq_if_in:NnF \g_cflist_loaded { #1 } {
		\seq_if_in:NnF \g_cflist_pending { #1 } {
			\seq_gput_right:Nn \g_cflist_pending { #1 }
		}
	}
}

\NewDocumentCommand{\cfconsiderloaded}{ m }{
	\seq_gput_right:Nn \g_cflist_loaded {#1}
}

\NewDocumentCommand{\cfremove}{ m }
{
	\seq_gremove_all:Nn \g_cflist_pending { #1 }
}

\NewDocumentCommand{\cfload}{ o }
{
	\seq_if_empty:NTF \g_cflist_pending {\unskip} {
		(cf.\ \cref{\seq_use:Nn \g_cflist_pending {,}})\IfValueTF{#1}{#1~}{\unskip}
		\seq_gconcat:NNN \g_cflist_loaded \g_cflist_loaded \g_cflist_pending
		\seq_gclear:N \g_cflist_pending
	}
}

\NewDocumentCommand{\cfclear} {} {
	\seq_gclear:N \g_cflist_loaded
	\seq_gclear:N \g_cflist_pending
}

\NewDocumentCommand{\cfout}{ o }
{
	\seq_if_empty:NTF \g_cflist_pending {\unskip} {
		(cf.\ \cref{\seq_use:Nn \g_cflist_pending {,}})\IfValueTF{#1}{#1~}{\unskip}
		\seq_gclear:N \g_cflist_pending
	}
}

\NewDocumentCommand{\ifnocf} { m } {
	\seq_if_empty:NT \g_cflist_pending { #1 }
}

\prop_const_from_keyval:Nn \l__verbs
{
	show = shows ,
	imply = implies ,
	demonstrate = demonstrates ,
	prove = proves ,
	establish = establishes ,
	ensure = ensures ,
	assure = assures
}

\seq_const_from_clist:Nn \g_prove_mru {
	establish,
	demonstrate,
	prove,
	show,
	imply,
	ensure
}

\tl_new:N \g_wordtmp
\seq_new:N \l_mytmps

\cs_generate_variant:Nn \str_if_in:nnTF { nVTF }

\NewDocumentCommand{\prove}{ o }{
	\IfValueTF{#1}{
		\seq_clear:N \l_mytmps
		\seq_map_inline:Nn \g_prove_mru {
			\str_if_eq:nnTF {##1} {ensure} {
				\str_set:Nn \l_temps {n}
			} {
				\str_set:Nx \l_temps {\str_head_ignore_spaces:n {##1}}
			}
			\str_if_in:nVTF {#1} \l_temps {
				\seq_put_right:Nn \l_mytmps {##1}
			} { }
		}
		\seq_get_right:NN \l_mytmps \g_wordtmp
	} {
		\seq_get_right:NN \g_prove_mru \g_wordtmp
	}
	\tl_use:N \g_wordtmp
	\seq_gput_left:NV \g_prove_mru \g_wordtmp
	\seq_gremove_duplicates:N \g_prove_mru
	\IfValueF{#1}{~}
}

\NewDocumentCommand{\proves}{ o }{
	\IfValueTF{#1}{
		\seq_clear:N \l_mytmps
		\seq_map_inline:Nn \g_prove_mru {
			\str_if_eq:nnTF {##1} {ensure} {
				\str_set:Nn \l_temps {n}
			} {
				\str_set:Nx \l_temps {\str_head_ignore_spaces:n {##1}}
			}
			\str_if_in:nVTF {#1} \l_temps {
				\seq_put_right:Nn \l_mytmps {##1}
			} { }
		}
		\seq_get_right:NN \l_mytmps \g_wordtmp
	} {
		\seq_get_right:NN \g_prove_mru \g_wordtmp
	}
	\str_set:NV \l_tmpa_str \g_wordtmp
	\prop_get:NVN \l__verbs \l_tmpa_str \l_tmpa_tl
	\tl_use:N \l_tmpa_tl
	\seq_gput_left:NV \g_prove_mru \g_wordtmp
	\seq_gremove_duplicates:N \g_prove_mru
	\IfValueF{#1}{~}
}

\ExplSyntaxOff

\NewDocumentEnvironment{cproof}{m}
{\begin{proof}[Proof of \cref{#1}]}%
	{\noindent The proof of \cref{#1} is thus complete.
\end{proof}}

\NewDocumentEnvironment{cproof2}{m}
{\begin{proof}[Proof of \cref{#1}]}%
	{\noindent This completes the proof of \cref{#1}.
\end{proof}}

\title{Normalized gradient flow optimization in the\\ training of 
ReLU artificial neural networks}

\author{
	Simon Eberle$^{1,2,a}$,
	Arnulf Jentzen$^{2,3,b}$,
	Adrian Riekert$^{3,c}$,
	and
	Georg Weiss$^{4,d}$
	\bigskip
	\\
	\small{$^1$\,Basque Center for Applied Mathematics, Bilbao, Spain}
	\\
	\small{$^2$\,School of Data Science and Shenzhen Research Institute of Big Data,}
	\\[-0.13cm]
	\small{The Chinese University of Hong Kong, Shenzhen, China}
	\\
	\small{$^3$\,Applied Mathematics: Institute for Analysis and Numerics,}
	\\[-0.13cm]
	\small{University of M{\"u}nster, Germany}
	\\
	\small{$^4$\,University of Duisburg-Essen, Germany} 
	\\
	\small{$^a$\,\textit{seberle@bcamath.org}}
	\\[-0.13cm]
	\small{$^b$\,\textit{ajentzen@uni-muenster.de}, \textit{ajentzen@cuhk.edu.cn}}
	\\[-0.13cm]
	\small{$^c$\,\textit{ariekert@uni-muenster.de}}
	\\[-0.13cm]
	\small{$^d$\,\textit{georg.weiss@uni-due.de}}
}

\date{\today}

\begin{document}

\maketitle

\begin{abstract}
	The training of artificial neural networks (ANNs) is nowadays 
	a highly relevant algorithmic procedure 
	with many applications in science and industry. 
	Roughly speaking, ANNs can be regarded 
	as iterated compositions between 
	affine linear functions and certain fixed nonlinear functions, 
	which are usually multidimensional versions 
	of a one-dimensional so-called activation function. 
	The most popular choice of such a one-dimensional 
	activation function is the rectified linear unit (ReLU) 
	activation function which maps a real number to its positive part 
	$ \R \ni x \mapsto \max\{ x, 0 \} \in \R $. 
	In this article we propose and analyze a modified variant 
	of the standard training procedure of such ReLU ANNs 
	in the sense that we propose to restrict the negative gradient flow 
	dynamics to a large submanifold of the ANN parameter space, 
	which is a strict $ C^{ \infty } $-submanifold 
	of the entire ANN parameter space that 
	seems to enjoy better regularity properties than the entire 
	ANN parameter space but which is also sufficiently large and sufficiently high dimensional 
	so that it can represent
	all ANN realization functions 
	that can be represented through the entire ANN parameter space. 
	In the special situation of shallow ANNs with 
	just one-dimensional ANN layers we also prove 
	for every Lipschitz continuous target function that 
	every gradient flow trajectory on this large submanifold 
	of the ANN parameter space is globally bounded. 
	For the standard gradient flow on the entire ANN parameter space 
	with Lipschitz continuous target functions 
	it remains an open problem of research to prove or disprove 
	the global boundedness of gradient flow trajectories even 
	in the situation of shallow ANNs with just one-dimensional 
	ANN layers. 
\end{abstract}

\tableofcontents

\section{Introduction}

The training of deep artificial neural networks (ANNs) is nowadays 
a highly relevant technical procedure 
with many applications in science and industry. 
In the most simple form we can think of a deep ANN as 
a tuple of real numbers describing a function, 
the so-called realization function of the ANN, 
which consists of multiple compositions of 
affine linear functions and certain fixed nonlinear functions. 
To be more specific, the realization function of such an ANN 
with 
$ L \in \N \cap (1,\infty) = \{ 2, 3, 4, \dots \} $ affine linear transformations 
and layer dimensions $ \ell_0, \ell_1, \dots, \ell_L \in \N = \{ 1, 2, 3, \dots \} $
is given through an affine linear function from $ \R^{ \ell_0 } $
to $ \R^{ \ell_1 } $ (1st affine linear transformation), 
then a fixed nonlinear function 
from $ \R^{ \ell_1 } $ to $ \R^{ \ell_1 } $, 
then again an affine linear function from $ \R^{ \ell_1 } $ 
to $ \R^{ \ell_2 } $ (2nd affine linear transformation), 
then again a fixed nonlinear function 
from $ \R^{ \ell_2 } $ to $ \R^{ \ell_2 } $, 
$ \dots $, and, finally, an affine linear 
function from $ \R^{ \ell_{ L - 1 } } $
to $ \R^{ \ell_L } $ ($ L $-th affine linear 
transformation). 
There are thus $ \ell_0 \ell_1 + \ell_1 $ real numbers 
to describe the 1st affine linear transformation 
in the ANN 
(affine linear transformation
from $ \R^{ \ell_0 } $ to $ \R^{ \ell_1 } $),
there are thus $ \ell_1 \ell_2 + \ell_2 $ real numbers 
to describe the 2nd affine linear transformation 
in the ANN 
(affine linear transformation
from $ \R^{ \ell_1 } $ to $ \R^{ \ell_2 } $), 
$ \dots $, 
and there are thus 
$ \ell_{ L - 1 } \ell_L + \ell_L $ real numbers 
to describe the $ L $-th affine linear transformation 
(affine linear transformation 
from $ \R^{ \ell_{ L - 1 } } $ to $ \R^{ \ell_L } $). 
The overall number $ \fd \in \N $ of 
real ANN parameters thus satisfies 
\begin{equation}
\textstyle 
\fd = 
\sum_{ k = 1 }^L
( \ell_{ k - 1 } \ell_k + \ell_k )
=
\sum_{ k = 1 }^L
\ell_k ( \ell_{ k - 1 } + 1 )
.
\end{equation}
We also refer to \cref{figure_dnn_general} for a graphical illustration 
of the architecture of such an ANN. 

\def\layersep{2.5cm}
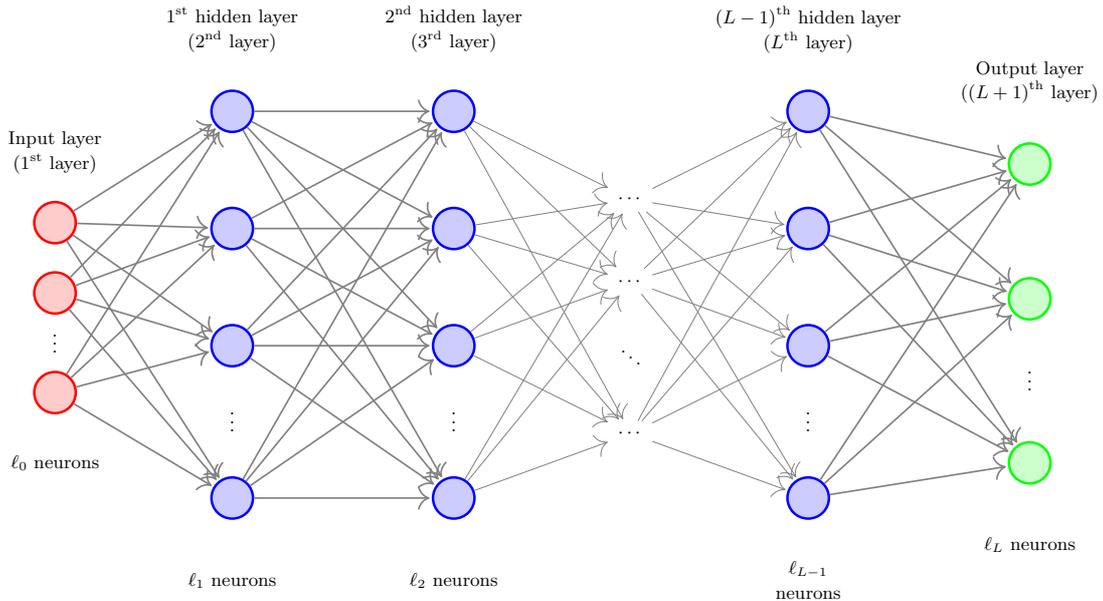
\begin{figure}[H]
	\centering
	\begin{adjustbox}{width=\textwidth}
		\begin{tikzpicture}[shorten >=1pt,->,draw=black!50, node distance=\layersep]
		\tikzstyle{every pin edge}=[<-,shorten <=1pt]
		\tikzstyle{input neuron}=[very thick, circle,draw=red, fill=red!20, minimum size=20pt,inner sep=0pt]
		\tikzstyle{output neuron}=[very thick, circle, draw=green,fill=green!20,minimum size=20pt,inner sep=0pt]
		\tikzstyle{hidden neuron}=[very thick, circle,draw=blue,fill=blue!20,minimum size=20pt,inner sep=0pt]
		\tikzstyle{annot} = [text width=9em, text centered]
		\tikzstyle{annot2} = [text width=4em, text centered]
		
		
		\node[input neuron] (I-1) at (-3,-1.4) {};
		\node[input neuron] (I-2) at (-3,-2.6) {};
		\node(I-dots) at (-3,-3.36) {\vdots};
		\node[input neuron] (I-3) at (-3,-4.3) {};
		
		\path[yshift = 1.5cm]
		node[hidden neuron](H0-1) at (0*\layersep, -1 cm) {};
		\path[yshift = 1.5cm]
		node[hidden neuron](H0-2) at (0*\layersep, -3 cm) {};
		\path[yshift = 1.5cm]
		node[hidden neuron](H0-3) at (0*\layersep, -5 cm) {};
		\path[yshift = 1.5cm]
		node(H0-dots) at (0*\layersep, -6.2 cm) {\vdots};
		\path[yshift = 1.5cm]
		node[hidden neuron](H0-4) at (0*\layersep, -7.6 cm) {};
		
		\path[yshift = 1.5cm]
		node[hidden neuron](H1-1) at (1.5*\layersep, -1 cm) {};
		\path[yshift = 1.5cm]
		node[hidden neuron](H1-2) at (1.5*\layersep, -3 cm) {};
		\path[yshift = 1.5cm]
		node[hidden neuron](H1-3) at (1.5*\layersep, -5 cm) {};
		\path[yshift = 1.5cm]
		node(H1-dots) at (1.5*\layersep, -6.2 cm) {\vdots};
		\path[yshift = 1.5cm]
		node[hidden neuron](H1-4) at (1.5*\layersep, -7.6 cm) {};
		
		\path[yshift = 0.5cm]
		node(Hdot-1) at (2.7*\layersep, -1.5 cm) {$\cdots$};
		\path[yshift = 0.5cm]
		node(Hdot-2) at (2.7*\layersep, -2.9 cm) {$\cdots$};
		\path[yshift = 0.5cm]
		node(Hdot-dots) at (2.7*\layersep, -4.1 cm) {$\ddots$};
		\path[yshift = 0.5cm]
		node(Hdot-3) at (2.7*\layersep, -5.5 cm) {$\cdots$};
		
		\path[yshift = 1.5cm]
		node[hidden neuron](H2-1) at (3.9*\layersep, -1 cm) {};
		\path[yshift = 1.5cm]
		node[hidden neuron](H2-2) at (3.9*\layersep, -3 cm) {};
		\path[yshift = 1.5cm]
		node[hidden neuron](H2-3) at (3.9*\layersep, -5 cm) {};
		\path[yshift = 1.5cm]
		node(H2-dots) at (3.9*\layersep, -6.2 cm) {\vdots};
		\path[yshift = 1.5cm]
		node[hidden neuron](H2-4) at (3.9*\layersep, -7.6 cm) {};

		\path[yshift = 1.5cm]
		node[output neuron](O-1) at (5.4*\layersep,-1.9 cm) {}; 
		\path[yshift = 1.5cm]
		node[output neuron](O-2) at (5.4*\layersep,-4.2 cm) {}; 
		\path[yshift = 1.5cm]
		node(O-dots) at (5.4*\layersep, -5.5 cm) {\vdots};
		\path[yshift = 1.5cm]
		node[output neuron](O-3) at (5.4*\layersep,-7 cm) {};
		\foreach \source in {1,2,3}
		\foreach \dest in {1,2,3,4}
		\path[-{>[length=2mm, width=4mm]}, line width = 0.8] (I-\source) edge (H0-\dest);

		\foreach \source in {1,2,3,4}
		\foreach \dest in {1,2,3,4}
		\path[-{>[length=2mm, width=4mm]}, line width = 0.8] (H0-\source) edge (H1-\dest);

		\foreach \source in {1,2,3,4}
		\foreach \dest in {1,2,3}
		\draw[-{>[length=2mm, width=4mm]}, path fading=east] (H1-\source) -- (Hdot-\dest);
		
		\foreach \source in {1,2,3}
		\foreach \dest in {1,2,3,4}
		\draw[-{>[length=2mm, width=4mm]}, path fading=west] (Hdot-\source) -- (H2-\dest);
		
		\foreach \source in {1,2,3,4}
		\foreach \dest in {1,2,3}
		\path[-{>[length=2mm, width=4mm]}, line width = 0.8] (H2-\source) edge (O-\dest);

		\node[annot,above of=H0-1, node distance=1.4cm, align=center] (hl) {$1^{\text{st}}$ hidden layer\\($2^{\text{nd}}$ layer)};
		\node[annot,above of=H1-1, node distance=1.4cm, align=center] (hl) {$2^{\text{nd}}$ hidden layer\\($3^{\text{rd}}$ layer)};
		\node[annot,above of=H2-1, node distance=1.4cm, align=center] (hl2) {${(L - 1)}^{\text{th}}$ hidden layer\\($L^{\text{th}}$ layer)};
		\node[annot,above of=I-1, node distance=1.2cm, align=center] {Input layer\\ ($1^{\text{st}}$ layer)};
		\node[annot,above of=O-1, node distance=1.4cm, align=center] {Output layer\\(${(L + 1)}^{\text{th}}$ layer)};
		
		\node[annot2,below of=H0-4, node distance=1.4cm, align=center] (sl) {$\ell_1$ neurons};
		\node[annot2,below of=H1-4, node distance=1.4cm, align=center] (sl) {$\ell_2$ neurons};
		\node[annot2,below of=H2-4, node distance=1.4cm, align=center] (sl2) {$\ell_{L - 1}$ neurons};
		\node[annot2,below of=I-3, node distance=1.2cm, align=center] {$\ell_0$ neurons};
		\node[annot2,below of=O-3, node distance=1.4cm, align=center] {$\ell_L$ neurons};
		\end{tikzpicture}
	\end{adjustbox}
	\caption{Graphical illustration for the architecture of an 
	ANN with $ L \in \N \cap (1,\infty) $ affine linear transformations, with 
	$ \ell_0 \in \N $ neurons on the input layer, with $ \ell_1 $ neurons on the $ 1 $\textsuperscript{st} hidden layer, 
	with $ \ell_2 $ neurons on the $ 2 $\textsuperscript{nd} hidden layer, 
	$ \dots $,
	with $ \ell_{ L - 1 } $ neurons on the $ ( L - 1 ) $\textsuperscript{th} hidden layer, 
	and 
	with $ \ell_L $ neurons on the output layer.}
	\label{figure_dnn_general}
\end{figure}

The nonlinear functions in between the affine linear 
transformation are usually multi-dimensional versions 
of a fixed one-dimensional function 
$ a \colon \R \to \R $ 
in the sense that the 
nonlinear function after the $ k $-th 
affine linear transformation 
with $ k \in \{ 1, 2, \dots, L - 1 \} $ 
is the function from 
$ \R^{ \ell_k } $
to
$ \R^{ \ell_k } $ 
given by
\begin{equation}
\label{eq:multidimensional_version}
\R^{ \ell_k } \ni ( x_1, \dots, x_{ \ell_k } ) 
\mapsto 
(
a( x_1 ), \dots, a( x_{ \ell_k } )
)
\in \R^{ \ell_k }
\end{equation}
and the one-dimensional function 
$ a \colon \R \to \R $ 
is then referred to as \emph{activation function} 
of the considered ANN. 
In numerical simulations maybe 
the most popular choice for the 
activation function 
$ a \colon \R \to \R $ 
in \cref{eq:multidimensional_version} 
is the ReLU activation function 
which is given by 
\begin{equation}
\label{eq:ReLU_activation}
\R \ni x \mapsto \max\{ x, 0 \} \in \R 
.
\end{equation}
There are also very good analytical reasons 
why the ReLU activation function 
in \cref{eq:ReLU_activation} 
seems to be so popular 
in numerical simulations. 
More formally, in the case of the ReLU activation 
in \cref{eq:ReLU_activation} 
it has been proven 
(see \cite{JentzenRiekertLojaDNN})
for Lipschitz continuous target functions 
that there exist global minimum points 
in the risk landscape in the training 
of ANNs in the shallow situation 
$ 
( L, \ell_0, \ell_1, \ell_2 ) 
\in \{ 2 \} \times \{ 1 \} \times \N \times \{ 1 \} 
$ 
while for other smooth activation functions 
such as the standard logistic activation 
function $ \R \ni x \mapsto ( 1 + \exp( - x ) )^{ - 1 } \in \R $
the existence of global minimum points 
has been disproven 
(see \cite{GallonJentzen2022,PetersenRaslanVoigtlaender2020}) 
and the existence of global minimum points in the risk landscape, in turn, 
seems to be closely related to the 
boundedness of gradient descent (GD) trajectories; 
see \cite{GallonJentzen2022}.

Despite the common usage of the 
ReLU activation function 
in deep ANNs, 
it remains an open problem of research 
to rigorously prove (or disprove) 
the convergence of GD trajectories. 
This lack of theoretical understanding 
applies to the ReLU activation function 
but also to other activation function. 
While for other smooth activation function, 
the boundedness of GD trajectories is often not even 
expected (see \cite{GallonJentzen2022}), for the ReLU activation function 
it remains an open problem to prove (or disprove) 
the boundedness of GD trajectories in the training 
of ReLU ANNs. 
Another key difficulty in the mathematical analysis of 
the training process of ReLU ANNs is the fact that 
the ReLU activation in \cref{eq:ReLU_activation}
fails to be differentiable at $ 0 $ 
and this lack of differentiability of 
the activation function transfers to 
the risk function, which, in turn, makes it difficult 
to analyze time-discrete GD processes 
as analyses of such methods rely on local Lipschitz 
continuity properties of the gradient of the risk function 
(see, e.g., \cite{BertsekasTsitsiklis2000,JentzenRiekertLojaDNN,Nesterov2004}).

In this article we propose and analyze a modified 
variant of the standard training process of ReLU ANNs. 
More formally, 
in this work we modify the usual gradient flow dynamics 
in a way so that the gradient flow remains the entire 
training process on a large 
submanifold 
of the ANN parameter space $ \R^{ \fd } $. 
Specifically, in this work we consider a suitable 
$ ( \fd - \sum_{ k = 1 }^{ L - 1 } \ell_k ) $-dimensional 
$ C^{ \infty } $-submanifold of the $ \fd $-dimensional ANN parameter 
space $ \R^{ \fd } $
and modify the gradient flow dynamics in a way so that the modified 
gradient flow remains on this 
submanifold.

The advantages of this gradient descent dynamics on this 
$ ( \fd - \sum_{ k = 1 }^{ L - 1 } \ell_k ) $-dimensional 
$ C^{ \infty } $-submanifold 
of the ANN parameter space $ \R^{ \fd } $ 
are 
\begin{enumerate}[label=(\roman*)]
	\item 
	that the risk function seems to have 
	better differentiability properties than 
	on the whole ANN parameter space $ \R^{ \fd } $
	and 
	\item 
	that nearly all parameters on the 
	submanifold are 
	bounded and it thus seems to be easier to verify the boundedness 
	of gradient flow trajectories on this submanifold. 
\end{enumerate}
In particular, in the special 
shallow ANN situation 
$ ( L, \ell_0, \ell_1, \ell_2 ) = (2, 1, 1, 1) $
we rigorously prove 
for every Lipschitz continuous target function 
$ f $
the global boundedness of every
gradient flow trajectory; 
see \cref{theo:bounded:1neuron} in \cref{sec:one_neuron_analysis} below. 
For the standard gradient flow on the entire ANN parameter space $ \R^{ \fd } $
with Lipschitz continuous target functions 
it remains an open problem of research to prove or disprove 
the global boundedness of gradient flow trajectories even 
in the special shallow ANN situation 
shallow ANN situation 
$ ( L, \ell_0, \ell_1, \ell_2 ) = (2, 1, 1, 1) $.

Let us also add a few references which are more or less related to the approach proposed in this article. In a very vague sense the approach in this article is related to the famous \emph{batch normalization procedure} (see Ioffe \& Szegedy~\cite{IoffeSzegedy2015})
 in the sense that in the batch normalization approach the data processed through the different ANN layers are normalized in a certain sense while in this article not the data processed through the ANN layers but the ANN parameters itself are normalized in a suitable sense.

As mentioned above, the risk function along the modified gradient flow trajectory appears to have better smoothness properties than on the entire ANN parameter space.
This is due to the fact that, roughly speaking, the input parameters of each hidden neuron have constant non-zero norm along the entire trajectory.
 The fact that in the case of shallow ANNs with one hidden layer certain differentiability properties can be ensured if one assumes that the inner ANN parameters are bounded away from zero in a suitable sense has previously been observed in, e.g., Chizat \& Bach~\cite{ChizatBach2020}, Wojtowytsch~\cite{Wojtowytsch2020}, and \cite[Proposition 2.11]{JentzenRiekertFlow}.

The remainder of this article is organized as follows. 
In \cref{sec:normalized_gradient_flow} we describe 
the modified gradient flow optimization dynamics 
in the situation of general deep ANNs 
with an arbitrary large number $ L \in \N \cap (1,\infty) $ 
of affine linear transformations 
and arbitrary layer dimensions 
$ \ell_0, \ell_1, \dots, \ell_L \in \N $. 
In \cref{sec:one_neuron_analysis} 
we consider the special situation 
of shallow ANNs with one-dimensional 
layer dimensions in the sense that 
$ L = 2 $ 
and 
$ \ell_0 = \ell_1 = \ell_2 = 1 $ 
and prove in \cref{theo:bounded:1neuron}
in this special situation 
for every Lipschitz continuous target function 
that every GF trajectory 
is globally bounded.

\section{Normalized gradient flow optimization 
	in the training of deep ReLU artificial neural networks (ANNs)}
\label{sec:normalized_gradient_flow}

In this section we describe and study 
the modified gradient flow optimization dynamics 
in the situation of general deep ANNs 
with an arbitrary large number $ L \in \N \cap (1,\infty) $ 
of affine linear transformations 
and arbitrary layer dimensions 
$ \ell_0, \ell_1, \dots, \ell_L \in \N $.

\subsection{Gradient flow optimization 
	on submanifolds on the ANN parameter space}

In the following abstract result, \cref{lem:modified_gradient},
we introduce a modification of a standard gradient flow $(\Theta_t ) _ { t \in [0 , \tau ) }$ with $\frac{\d}{\d t } \Theta_t = \cG ( \Theta_t ) $ (see \cref{eq:gradient_manifold_dynamics})
with the property that certain quantities $\psi_k ( \Theta_t ) \in \R$, $k \in \cu{1, 2, \ldots, K }$, are time-invariant. 
Roughly speaking, at each time $t \in [0 , \tau )$ the derivative vector $\cG ( \Theta_t) \in \R^\fd$ is projected onto the tangent space to a certain submanifold of $\R^\fd$ on which all $\psi_k$ are constant.
 Intuitively, this causes the gradient flow to move only tangentially to the manifold and therefore the quantities $\psi_k ( \Theta_t )$ remain invariant.

\begin{lemma}[Gradient flow dynamics on submanifolds]
	\label{lem:modified_gradient}
	Let $ \cG \colon \R^{ \fd } \to \R^{ \fd } $ be measurable, 
	let $ K \in \N $, 
	for every $ k \in \{ 1, 2, \dots, K \} $ 
	let $ \psi_k \colon \R^{ \fd } \to \R $ 
	be continuously differentiable, 
	assume for all $ k, l \in \{ 1, 2, \dots, K \} $, 
	$ \theta \in \R^{ \fd } $
	with 
	$
	\min_{ m \in \{ 1, 2, \dots, K \} }
	\| ( \nabla \psi_m )( \theta ) \| > 0
	$
	and 
	$ k \neq l $ that 
	\begin{equation}
	\label{eq:scalar_product_orthogonality}
	\spro{ ( \nabla \psi_k )( \theta ), ( \nabla \psi_l )( \theta ) }
	= 0
	,
	\end{equation}
	let $\tau \in (0 , \infty ]$,
	and let $ \Theta \in C( [0 , \tau ), \R^{ \fd } ) $ 
	satisfy for all $ t \in [0 , \tau ) $ that 
	$
	\inf_{ s \in [0,t] } 
	\min_{ k \in \{ 1, 2, \dots, K \} } \allowbreak
	\| ( \nabla \psi_k )( \Theta_s ) \| 
	\allowbreak
	> 0
	$,
	$
	\int_0^t
	\| \cG( \Theta_s ) \|
	\, \d s
	< \infty
	$, 
	and 
	\begin{equation}
	\label{eq:gradient_manifold_dynamics}
	\Theta_t 
	=
	\Theta_0
	+
	\int\limits_0^t
	\rbr*{
	\cG( \Theta_s )
	-
	\ssuml_{ k = 1 }^K
	\| ( \nabla \psi_k )( \Theta_s ) \|^{ - 2 }
	\langle
	\cG( \Theta_s ), 
	( \nabla \psi_k )( \Theta_s )
	\rangle
	( \nabla \psi_k )( \Theta_s ) } \,
	\d s .
	\end{equation}
	Then 
	it holds for all 
	$ k \in \{ 1, 2, \dots, K \} $, 
	$ t \in [0, \tau ) $ 
	that
	$
	\psi_k( \Theta_t ) = \psi_k( \Theta_0 )
	$.
\end{lemma}
\begin{cproof}{lem:modified_gradient}
	\Nobs that \cref{eq:scalar_product_orthogonality} 
   \proves for all 
	$ k \in \{ 1, 2, \dots, K \} $, 
	$ \theta \in \R^{ \fd } $
	with 
	$
	\min_{ l \in \{ 1, 2, \dots, K \} }
	\| ( \nabla \psi_l )( \theta ) \| > 0
	$
	that 
	\begin{equation}
	\begin{split}
	&
	\textstyle 
	\spro*{
	( \nabla \psi_k )( \theta ) ,
	\cG( \theta )
	-
	\sum\limits_{ l = 1 }^K
	\| ( \nabla \psi_l )( \theta ) \|^{ - 2 }
	\langle
	\cG( \theta ), 
	( \nabla \psi_l )( \theta )
	\rangle
	( \nabla \psi_l )( \theta ) }
	\\ &
	=
	\textstyle 
	\bigl\langle 
	( \nabla \psi_k )( \theta ) ,
	\cG( \theta )
	\bigr\rangle 
	-
	\sum\limits_{ l = 1 }^K
	\Bigl\langle 
	( \nabla \psi_k )( \theta ) ,
	\| ( \nabla \psi_l )( \theta ) \|^{ - 2 }
	\langle
	\cG( \theta ), 
	( \nabla \psi_l )( \theta )
	\rangle
	( \nabla \psi_l )( \theta )
	\Bigr\rangle 
	\\ &
	=
	\textstyle 
	\bigl\langle 
	( \nabla \psi_k )( \theta ) ,
	\cG( \theta )
	\bigr\rangle 
	-
	\sum\limits_{ l = 1 }^K
	\Bigl[
	\| ( \nabla \psi_l )( \theta ) \|^{ - 2 }
	\langle
	\cG( \theta ), 
	( \nabla \psi_l )( \theta )
	\rangle
	\Bigr]
	\langle 
	( \nabla \psi_k )( \theta ) ,
	( \nabla \psi_l )( \theta )
	\rangle 
	\\ &
	=
	\textstyle 
	\bigl\langle 
	( \nabla \psi_k )( \theta ) ,
	\cG( \theta )
	\bigr\rangle 
	-
	\| ( \nabla \psi_k )( \theta ) \|^{ - 2 }
	\langle
	\cG( \theta ), 
	( \nabla \psi_k )( \theta )
	\rangle
	\langle 
	( \nabla \psi_k )( \theta ) ,
	( \nabla \psi_k )( \theta )
	\rangle 
	\\ &
	=
	\textstyle 
	\bigl\langle 
	( \nabla \psi_k )( \theta ) ,
	\cG( \theta )
	\bigr\rangle 
	-
	\langle
	\cG( \theta ), 
	( \nabla \psi_k )( \theta )
	\rangle
	= 0
	.
	\end{split}
	\end{equation}
	The generalized chain rule 
	and \cref{eq:gradient_manifold_dynamics} 
	hence \prove that for all 
	$ k \in \{ 1, 2, \dots, K \} $, 
	$ t \in [0,\tau ) $
	it holds that
	\begin{equation}
	\psi_k( \Theta_t )
	=
	\psi_k( \Theta_0 )
	.
	\end{equation}
\end{cproof}

\subsection{Descent property for modified gradient flows}

In this subsection we show in an abstract setting that the considered modified gradient flow still has a descent property in the sense that the value of the objective function $\cL  ( \Theta_t )$, $t \in [0 , \infty )$,
 is monotonically non-increasing in time.
Notice that we do not assume that the objective function $\cL \colon U \to \R$ is continuously differentiable. 
Instead, we only assume that $\cL$ can be approximated by differentiable functions $\fL_r \in C^1 ( U , \R )$, $r \in \N$, in a suitable sense (see below \cref{eq:modified_gradient_dynamics} in \cref{prop:descent0}).
This will be important when applying our results to the risk functions occurring in the training of ANNs with the non-differentiable ReLU activation.
For the proof of \cref{prop:descent0} we will apply the generalized chain rule from Cheridito et al.~\cite[Lemma 3.3]{CheriditoJentzenRiekert2022}.

\begin{prop}[Energy dynamics for modified gradient flows]
	\label{prop:descent0}
	Let $ \fd, K \in \N $, 
	for every $ k \in \{ 1, 2, \dots, K \} $
	let 
	$ \phi_k \colon \R^{ \fd } \to \R $ 
	and 
	$ \psi_k \colon \R^{ \fd } \to \R^{ \fd } $ 
	be locally bounded and measurable, 
	let $ U \subseteq \R^{ \fd } $ 
	be open, 
	let 
	$ \gamma \in C( [0,\infty), [0,\infty) ) $, 
	$ \cL \in C( U, \R ) $, let $\cG \colon U \to \R^\fd$ be locally bounded and measurable,
	let
	$ \Theta \in C( [0,\infty), U ) $ 
	satisfy for all $ t \in [0,\infty) $ that 
	\begin{equation}
	\label{eq:modified_gradient_dynamics}
	\Theta_t
	=
	\Theta_0
	-
	\int_0^t
	\textstyle
	\gamma(s)
	\,
	\Bigl(
	\cG ( \Theta_s )
	+
	\sum\limits_{ k = 1 }^K
	\phi_k( \Theta_s )
	\langle \psi_k( \Theta_s ), \cG ( \Theta_s ) \rangle 
	\psi_k( \Theta_s )
	\Bigr)
	\, \d s,
	\end{equation}
	and assume that there exist $\fL_r \in C^1 ( U , \R^\fd )$, $r \in \N$,
	which satisfy for all compact $K \subseteq U$ that
	$\sup_{r \in \N} \sup_{\theta \in K } \norm{ \nabla \fL_r ( \theta ) } < \infty $
	and which satisfy for all $\theta \in U$ that $\lim_{r \to \infty } \fL_r ( \theta ) = \cL ( \theta ) $ and $\lim_{r \to \infty} \nabla \fL_r ( \theta ) = \cG ( \theta ) $.
	Then it holds for all $t \in [0 , \infty )$ that
	\begin{equation}
	\label{eq:gradient_flow_accelerated}
	\cL( \Theta_t )
	=
	\cL( \Theta_0 )
	-
	\int_0^t
	\gamma(s)
	\,
	\Bigl(
	\| \cG ( \Theta_s ) \|^2
	+
	\textstyle 
	\sum\limits_{ k = 1 }^K
	\phi_k( \Theta_s )
	|
	\langle \psi_k( \Theta_s ), \cG ( \Theta_s ) \rangle 
	|^2
	\Bigr)
	\, \d s
	.
	\end{equation}
\end{prop}
\begin{cproof}{prop:descent0}
	\Nobs that the assumption that $\phi_k$, $k \in \cu{1, 2, \ldots, K }$, $\psi_k$, $k \in \cu{1, 2, \ldots, K }$, and $\cG$ are locally bounded and measurable and the fact that $\Theta$ is continuous \prove for all $t \in [0 , \infty )$
	that
	$
	[0, t ] \ni s \mapsto \gamma(s) 
	( \cG ( \Theta_s )
	+
	\sum_{ k = 1 }^K
	\phi_k( \Theta_s )
	\langle \psi_k( \Theta_s ), \cG ( \Theta_s ) \rangle 
	\psi_k( \Theta_s ) ) \in \R^\fd
	$
	is bounded and measurable. 
	Combining this with \cref{eq:modified_gradient_dynamics} and the generalized chain rule (cf., e.g., Cheridito et al.~\cite[Lemma 3.3]{CheriditoJentzenRiekert2022}) \proves 
	 for all $t \in [0 , \infty )$, $r \in \N$ that
	\begin{equation}
	\begin{split}
	 \fL_r ( \Theta_t ) &= \fL_r ( \Theta_0 ) 
 - \int_0^t \gamma( s ) \Bigl(  \spro{\cG ( \Theta_s ) , \nabla \fL_r ( \Theta_s ) } \\
 & \quad  + \ssuml_{k=1}^K \phi_k ( \Theta_s ) \spro{ \psi_k( \Theta_s ), \cG ( \Theta_s ) } \spro{ \psi_k ( \Theta_s ) , \nabla \fL_r ( \Theta_s ) } \Bigr) \, \d s.
	\end{split}
	\end{equation}
	In addition, \nobs that the fact that $\Theta$ is continuous \proves 
	 for every $t \in [0 , \infty )$ that $\cu{ \Theta_s \colon s \in [0 , t ] } \subseteq U$ is compact.
	Combining this with the assumption that for all compact $K \subseteq U$ it holds that
	$\sup_{r \in \N} \sup_{\theta \in K } \norm{ \nabla \fL_r ( \theta ) } < \infty $,
	the assumption that for all $\theta \in K$ it holds that $\lim_{r \to \infty } \fL_r ( \theta ) = \cL ( \theta ) $ and $\lim_{r \to \infty} \nabla \fL_r ( \theta ) = \cG ( \theta ) $,
	and the dominated convergence theorem \proves \cref{eq:gradient_flow_accelerated}.
\end{cproof}

\begin{cor}[Energy dynamics for modified gradient flows]
	\label{cor:descent2}
	Let $ \fd, K \in \N $, 
	for every $ k \in \{ 1, 2, \dots, K \} $
	let 
	$ \psi_k \colon \R^{ \fd } \to \R^{ \fd } $ 
	be locally bounded and measurable, 
	let $ U \subseteq \R^{ \fd } $ 
	be open, 
	let 
	$ \gamma \in C( [0,\infty), [0,\infty) ) $, 
	$ \cL \in C( U, \R ) $, let $\cG \colon U \to \R^\fd$ be locally bounded and measurable,
	let
	$ \Theta \in C( [0,\infty), U ) $ 
	satisfy for all $ t \in [0,\infty) $ that 
	\begin{equation}
	\label{eq:modified_gradient_dynamics2}
	\Theta_t
	=
	\Theta_0
	-
	\int_0^t
	\textstyle
	\gamma(s) \, 
	\Bigl(
	\cG ( \Theta_s )
	+
	\sum\limits_{ k = 1 }^K
	\langle \psi_k( \Theta_s ), \cG( \Theta_s ) \rangle 
	\psi_k( \Theta_s )
	\Bigr)
	\, \d s,
	\end{equation}
	and assume that there exist $\fL_r \in C^1 ( U , \R^\fd )$, $r \in \N$,
	which satisfy for all compact $K \subseteq U$ that
	$\sup_{r \in \N} \sup_{\theta \in K } \norm{ \nabla \fL_r ( \theta ) } < \infty $
	and which satisfy for all $\theta \in U$ that $\lim_{r \to \infty } \fL_r ( \theta ) = \cL ( \theta ) $ and $\lim_{r \to \infty} \nabla \fL_r ( \theta ) = \cG ( \theta ) $.
	Then it holds for all $t \in [0 , \infty )$ that
	\begin{equation}
	\label{eq:gradient_flow_accelerated2}
	\cL( \Theta_t )
	=
	\cL( \Theta_0 )
	-
	\int_0^t
	\gamma(s) \,
	\Bigl(
	\| \cG ( \Theta_s ) \|^2
	+
	\textstyle 
	\sum\limits_{ k = 1 }^K
	|
	\langle \psi_k( \Theta_s ), \cG ( \Theta_s ) \rangle 
	|^2
	\Bigr)
	\, \d s
	.
	\end{equation}
\end{cor}
\begin{cproof}{cor:descent2}
	\Nobs that 
	\cref{prop:descent0}
	\proves[e] \cref{eq:gradient_flow_accelerated2}. 
\end{cproof}

In the next result we apply the more general \cref{prop:descent0} to the modified gradient flow from \cref{lem:modified_gradient}. 
Using Parseval's identity for the orthogonal gradient vectors $\nabla \psi_k ( \Theta_t ) \in \R^\fd$, $k \in \cu{1, 2, \ldots, K }$,
we establish that the value $\cL ( \Theta_t )$ is non-increasing in the time variable $t$.

\begin{cor}[Gradient flow dynamics on submanifolds]
	\label{cor:descent}
	Let $ \fd, K \in \N $, 
	for every $ k \in \{ 1, 2, \dots, K \} $ 
	let $ \psi_k \colon \R^{ \fd } \to \R $ 
	be continuously differentiable, 
	let $ U \subseteq \R^{ \fd } $ satisfy 
	\begin{equation}
	\textstyle 
	U = \{ 
	\theta \in \R^{ \fd } \colon
	\min_{ k \in \{ 1, 2, \dots, K \} }
	\| ( \nabla \psi_k )( \theta ) \| > 0
	\} 
	,
	\end{equation}
	let 
	$ \cL \in C( U, \R ) $, let $\cG \colon U \to \R^\fd$ be locally bounded and measurable,
	assume for all 
	$ \theta \in U $, 
	$ k, l \in \{ 1, 2, \dots, K \} $ 
	with $ k \neq l $ that 
	\begin{equation}
	\label{eq:scalar_product_orthogonality2}
	\spro{ ( \nabla \psi_k )( \theta ), ( \nabla \psi_l )( \theta ) }
	= 0	,
	\end{equation}
	let 
	$ \gamma \in C( [0,\infty), [0,\infty) ) $, 
	$ \Theta \in C( [0,\infty), U ) $ 
	satisfy for all $ t \in [0,\infty) $ that 
	\begin{equation}
	\label{eq:gradient_manifold_dynamics2}
	\textstyle 
	\Theta_t 
	=
	\Theta_0
	-
	\int\limits_0^t
	\gamma(s)
	\rbr*{
	\cG ( \Theta_s )
	-
	\sum\limits_{ k = 1 }^K
	\| ( \nabla \psi_k )( \Theta_s ) \|^{ - 2 }
	\langle
	\cG ( \Theta_s ), 
	( \nabla \psi_k )( \Theta_s )
	\rangle
	( \nabla \psi_k )( \Theta_s ) } \,
	\d s,
	\end{equation}
	and assume that there exist $\fL_r \in C^1 ( U , \R^\fd )$, $r \in \N$,
	which satisfy for all compact $K \subseteq U$ that
	$\sup_{r \in \N} \sup_{\theta \in K } \norm{ \nabla \fL_r ( \theta ) } < \infty $
	and which satisfy for all $\theta \in U$ that $\lim_{r \to \infty } \fL_r ( \theta ) = \cL ( \theta ) $ and $\lim_{r \to \infty} \nabla \fL_r ( \theta ) = \cG ( \theta ) $.
	Then 
	\begin{enumerate}[label=(\roman*)]
		\item 
		it holds for all 
		$ k \in \{ 1, 2, \dots, K \} $, 
		$ t \in [0,\infty) $ 
		that
		$
		\psi_k( \Theta_t ) = \psi_k( \Theta_0 )
		$
		and 
		\item 
		it holds for all $ s, t \in [0,\infty) $ 
		with $ s \leq t $ that 
		\begin{equation}
		\begin{split}
		&
		\cL( \Theta_t ) 
		=
		\cL( \Theta_s ) \\ &
		- \int_s^t
		\gamma ( u )
		\rbr[\Big]{
		\norm*{
		\cG ( \Theta_u ) } ^2
		-
		\ssuml_{ k = 1 }^K
		\| ( \nabla \psi_k )( \Theta_u ) \|^{ - 2 }
		|
		\langle
		\cG ( \Theta_u ), 
		( \nabla \psi_k )( \Theta_u )
		\rangle
		|^2 }
		\, \d u
		\\ &
		\textstyle 
		=
		\cL( \Theta_s )
		\\ &
		\textstyle 
		-
		\displaystyle\int_s^t
		\gamma ( u )
		\Bigl\| 
		\cG ( \Theta_u )
		-
		\ssuml_{ k = 1 }^K 
		\| ( \nabla \psi_k ) ( \Theta_u ) \|^{ - 2 }
		\langle
		\cG ( \Theta_s ), 
		( \nabla \psi_k ) ( \Theta_u )
		\rangle
		( \nabla \psi_k ) ( \Theta_u )
		\Bigr\|^2
		\d u 
		\\ &
		\leq 
		\cL( \Theta_s )
		.
		\end{split}
		\end{equation}
	\end{enumerate}
\end{cor}
\begin{cproof}{cor:descent}
	Throughout this proof 
	for every $ t \in [0,\infty) $
	let 
	$ P_t \colon \R^{ \fd } \to \R^{ \fd } $ 
	satisfy for all 
	$ v \in \R^{ \fd } $
	that 
	\begin{equation}
	P_t( v )
	=
	\sum_{ k = 1 }^K
	\|
	( \nabla \psi_k )( \Theta_t )
	\|^{ - 2 }
	\langle 
	( \nabla \psi_k )( \Theta_t ), v 
	\rangle 
	( \nabla \psi_k )( \Theta_t )
	.
	\end{equation}
	\Nobs that Parseval's identity 
	\proves that 
	for all $ t \in [0,\infty) $
	it holds that
	\begin{equation}
	\| P_t( v ) \|^2
	=
	\sum_{
		k = 1
	}^K
	\abs*{ \spro*{
	\frac{ ( \nabla \psi_k )( \Theta_t ) }{ \| ( \nabla \psi_k )( \Theta_t ) \| } , v } } ^2
	\leq 
	\| v \|^2
	.
	\end{equation}
	In addition,
	\cref{lem:modified_gradient}
	and 
	\cref{prop:descent0}
	\prove for all $s, t \in [0 , \infty )$ with $s \le t$ that 
	\begin{equation}
	\begin{split}
	&
	\cL( \Theta_t ) 
	=
	\cL( \Theta_s )
	\\ &
	-
	\int_s^t
	\gamma ( u )
	\rbr[\Big]{ \norm*{
	\cG ( \Theta_u ) } ^2
	-
	\ssuml_{ k = 1 }^K
	\| ( \nabla \psi_k )( \Theta_u ) \|^{ - 2 }
	|
	\langle
 \cG ( \Theta_u ), 
	( \nabla \psi_k )( \Theta_u )
	\rangle
	|^2 } 
	\, \d u
	\\ & =
	\cL( \Theta_s )
	-
	\int_s^t
	\gamma ( u )
	\rbr[\Big]{
	\norm*{
	\cG ( \Theta_u ) } ^2
	-
	\bigl\|
	P_u \bigl(
	\cG ( \Theta_u )
	\bigr)
	\bigr\|^2 } 
	\, \d u
	\\ & =
	\cL( \Theta_s )
	-
	\int_s^t
	\gamma( u )
	\bigl\| 
	\cG ( \Theta_u )
	-
	P_u \bigl(
	\cG ( \Theta_u )
	\bigr)
	\bigr\|^2
	\, \d u
	.
	\end{split}
	\end{equation}
\end{cproof}

Next, in \cref{cor:descent:c1} we specialize the above results to the case of a continuously differentiable objective function $\cL \in C^1 ( U , \R )$.

\begin{cor}[Gradient flow dynamics on submanifolds, differentiable case]
	\label{cor:descent:c1}
	Let $ \fd, K \in \N $, 
	for every $ k \in \{ 1, 2, \dots, K \} $ 
	let $ \psi_k \colon \R^{ \fd } \to \R $ 
	be continuously differentiable, 
	let $ U \subseteq \R^{ \fd } $ satisfy 
	\begin{equation}
	\textstyle 
	U = \{ 
	\theta \in \R^{ \fd } \colon
	\min_{ k \in \{ 1, 2, \dots, K \} }
	\| ( \nabla \psi_k )( \theta ) \| > 0
	\} 
	,
	\end{equation}
	let 
	$ \cL \in C^1 ( U, \R ) $,
	assume for all 
	$ \theta \in U $, 
	$ k, l \in \{ 1, 2, \dots, K \} $ 
	with $ k \neq l $ that 
	\begin{equation}
	\label{eq:scalar_product_orthogonality:c1}
	\spro{ ( \nabla \psi_k )( \theta ), ( \nabla \psi_l )( \theta ) }
	= 0 ,
	\end{equation}
	and let 
	$ \gamma \in C( [0,\infty), [0,\infty) ) $, 
	$ \Theta \in C( [0,\infty), U ) $ 
	satisfy for all $ t \in [0,\infty) $ that 
	\begin{equation}
	\label{eq:gradient_manifold_dynamics:c1}
	\textstyle 
	\Theta_t 
	=
	\Theta_0
	-
	\int\limits_0^t
	\gamma(s)
	\rbr*{
	( \nabla \cL ) ( \Theta_s )
	-
	\sum\limits_{ k = 1 }^K
	\| ( \nabla \psi_k )( \Theta_s ) \|^{ - 2 }
	\langle
	( \nabla \cL ) ( \Theta_s ), 
	( \nabla \psi_k )( \Theta_s )
	\rangle
	( \nabla \psi_k )( \Theta_s ) } \,
	\d s.
	\end{equation}
	Then 
	\begin{enumerate}[label=(\roman*)]
		\item 
		it holds for all 
		$ k \in \{ 1, 2, \dots, K \} $, 
		$ t \in [0,\infty) $ 
		that
		$
		\psi_k( \Theta_t ) = \psi_k( \Theta_0 )
		$
		and 
		\item 
		it holds for all $ s, t \in [0,\infty) $ 
		with $ s \leq t $ that 
		\begin{equation}
		\begin{split}
		&
		\cL( \Theta_t ) 
		=
		\cL( \Theta_s ) \\ &
		-
		\int_s^t
		\gamma ( u )
		\rbr[\Big]{
		\norm*{
		( \nabla \cL ) ( \Theta_s ) } ^2
		-
		\textstyle 
		\sum\limits_{ k = 1 }^K
		\| ( \nabla \psi_k )( \Theta_u ) \|^{ - 2 }
		|
		\langle
		( \nabla \cL )  ( \Theta_u ), 
		( \nabla \psi_k )( \Theta_u )
		\rangle
		|^2 } 
		\, \d u
		\\ &
		\textstyle 
		=
		\cL( \Theta_s )
		\\ &
		-
		\int_s^t
		\gamma ( u )
		\Bigl\| 
		( \nabla \cL ) ( \Theta_s )
		-
		\ssuml_{ k = 1 }^K 
		\| ( \nabla \psi_k ) ( \Theta_u ) \|^{ - 2 }
		\langle
		( \nabla \cL ) ( \Theta_u ), 
		( \nabla \psi_k ) ( \Theta_u )
		\rangle
		( \nabla \psi_k ) ( \Theta_u )
		\Bigr\|^2
		\d u
		\\ &
		\leq 
		\cL( \Theta_s )
		.
		\end{split}
		\end{equation}
	\end{enumerate}
\end{cor}

\begin{cproof}{cor:descent:c1}
	This is a special case of \cref{cor:descent}
	(applied with $\cL \with \cL$, $\cG \with \nabla \cL $, $(\fL_r)_{r \in \N } \with (\cL ) _{r \in \N }$ in the notation of \cref{cor:descent}).
\end{cproof}

In the final result of this subsection, \cref{cor:descent_B},
we show a modified version of \cref{cor:descent:c1}. More specifically, we prove that the time-dependent factors $\gamma(t) \in [0 , \infty )$, $t \in [0 , \infty )$, can be chosen in such a way that the value $\cL ( \Theta_t )$, $t \in [0 , \infty )$, decreases at the same rate as for the standard gradient flow; see \cref{cor:descent:c1:eq:claim} below.

\begin{cor}[Gradient flow dynamics on submanifolds]
	\label{cor:descent_B}
	Let $ \fd, K \in \N $, 
	for every $ k \in \{ 1, 2, \dots, K \} $ 
	let $ \psi_k \colon \R^{ \fd } \to \R $ 
	be continuously differentiable, 
	let $ U \subseteq \R^{ \fd } $ satisfy 
	\begin{equation}
	\textstyle 
	U = \{ 
	\theta \in \R^{ \fd } \colon
	\min_{ k \in \{ 1, 2, \dots, K \} }
	\| ( \nabla \psi_k )( \theta ) \| > 0
	\} 
	,
	\end{equation}
	assume for all 
	$ \theta \in U $, 
	$ k, l \in \{ 1, 2, \dots, K \} $ 
	with $ k \neq l $ that 
	\begin{equation}
	\label{eq:scalar_product_orthogonality2_B}
	\spro{ ( \nabla \psi_k )( \theta ), ( \nabla \psi_l )( \theta ) }
	= 0
	,
	\end{equation}
	and 
	let 
	$ \gamma \in C( [0,\infty), [0,\infty) ) $, 
	$ \cL \in C^1( U, \R ) $, 
	$ \Theta \in C( [0,\infty), U ) $ 
	satisfy for all $ t \in [0,\infty) $ that 
	\begin{equation} 
	\textstyle 
	\gamma(t)
	\biggl\| 
	( \nabla \cL )( \Theta_t )
	-
	\sum\limits_{ k = 1 }^K 
	\frac{ 
		\langle
		( \nabla \cL )( \Theta_t ), 
		( \nabla \psi_k )( \Theta_t )
		\rangle
		( \nabla \psi_k )( \Theta_t )
	}{
		\| ( \nabla \psi_k )( \Theta_t ) \|^2
	}
	\biggr\|^2
	=
	\| ( \nabla \cL )( \Theta_t ) \|^2
	\end{equation}
	and 
	\begin{equation}
	\label{eq:gradient_manifold_dynamics2_B}
	\textstyle 
	\Theta_t 
	=
	\Theta_0
	-
	\int\limits_0^t
	\gamma(s)
	\rbr*{
	( \nabla \cL )( \Theta_s )
	-
	\sum\limits_{ k = 1 }^K
	\| ( \nabla \psi_k )( \Theta_s ) \|^{ - 2 }
	\langle
	( \nabla \cL )( \Theta_s ), 
	( \nabla \psi_k )( \Theta_s )
	\rangle
	( \nabla \psi_k )( \Theta_s ) } \,
	\d s
	.
	\end{equation}
	Then 
	\begin{enumerate}[label=(\roman*)]
		\item 
		it holds for all 
		$ k \in \{ 1, 2, \dots, K \} $, 
		$ t \in [0,\infty) $ 
		that
		$
		\psi_k( \Theta_t ) = \psi_k( \Theta_0 )
		$
		and 
		\item \label{cor:descent:c1:eq:claim}
		it holds for all $ s, t \in [0,\infty) $ 
		with $ s \leq t $ that 
		\begin{equation}
		\begin{split}
		\cL( \Theta_t ) 
		=
		\cL( \Theta_s )
		-
		\int_s^t
		\norm*{
		( \nabla \cL )( \Theta_u ) } ^2
		\d u
		.
		\end{split}
		\end{equation}
	\end{enumerate}
\end{cor}
\begin{cproof}{cor:descent_B}
	This is a special case of \cref{cor:descent:c1}. 
\end{cproof}

\subsection{Normalized gradient descent  
	in the training of deep ReLU ANNs}

In the following we introduce our notation for deep ANNs with ReLU activation.
\cref{setting:dnn_normalized} below is inspired by Hutzenthaler et al.~\cite[Setting 2.1]{HutzenthalerJentzenPohlRiekertScarpa2021}.

In \cref{setting:dnn_normalized} we first introduce
the depth $ L \in \N \cap (1,\infty) $ 
of the considered ANN,
the layer dimensions 
$ \ell_0, \ell_1, \dots, \ell_L \in \N $,
 the continously differentiable approximations $\Rect_r \in C ^1 ( \R , \R )$, $r \in \N$, for the ReLU activation function $\Rect_{ \infty } ( x ) = \max \cu{x , 0 } $,
the unnormalized probability distribution $ \mu \colon \cB( [ a, b ]^{ \ell_0 } ) \to [0, \infty] $ of the input data,
and the measurable target function $f \colon [a,b]^{ \ell_0 } \to \R^{ \ell_L }  $.
Note that in the definition of the ANN realization functions in \cref{eq:def_NN_realization_normalized} we subtract in the last layer (the case $k+1 = L$) the average value of the output of the previous layer with respect to the input distribution $\mu$.
In \cref{eq:def_risk_function_normalized} we introduce the risk functions $\cL_r \colon \R^\fd \to \R$, $r \in \N \cup \cu{ \infty }$, and we define the generalized gradient $\cG \colon \R^\fd \to \R^\fd$ as the pointwise limit of the approximate gradients $\nabla \cL_r \colon \R^\fd \to \R^\fd$ as $r \to \infty $.

In \cref{eq:setting:def:psik} we inductively define the layer-wise rescaling operations $\Psi_k \colon \R^\fd \to \R^\fd$, $k \in \cu{0 , 1, \ldots, L - 1 }$,
 which have the property that certain sub-vectors of the parameter vector $\Psi_{L - 1 } ( \theta )$ are modified in order to have norm $1$ without changing the realization function; see \cref{prop:rescaling:properties} below for details.
 Finally, we define the modified gradient flow process $\Theta \colon [0 , \infty ) \times \Omega \to \R^\fd$ with random initialization $\xi \colon \Omega \to \R^\fd$
 and the modified gradient descent process $\varTheta \colon \N_0 \times \Omega \to \R^\fd$.

\begin{setting}
	\label{setting:dnn_normalized}
	Let 
	$ a \in \R $, 
	$ b \in (a,\infty) $, 
	$ ( \ell_k )_{ k \in \N_0 } \subseteq \N $, 
	$ L, \fd \in \N \backslash \{ 1 \} $
	satisfy 
	$
	\mathfrak{d} = \sum_{ k = 1 }^L \ell_k ( \ell_{ k - 1 } + 1 ) 
	$, 
	for every 
	$ 
	\theta = ( \theta_1, \dots, \theta_{ \fd } ) \in \R^{ \fd } 
	$
	let 
	$ 
	\fw^{ k, \theta } = 
	( \fw^{ k, \theta }_{ i, j } )_{ 
		(i,j) \in \{ 1, \ldots, \ell_k \} \times \{ 1, \ldots, \ell_{ k - 1 } \} 
	}
	\in \R^{ \ell_k \times \ell_{ k - 1 } }
	$, 
	$
	k \in \N 
	$, 
	and 
	$
	\fb^{ k, \theta } 
	= 
	( \fb^{ k, \theta }_1, \dots, \fb^{ k, \theta }_{ \ell_k} )
	\in \R^{ \ell_k } 
	$,
	$ k \in \N $, 
	satisfy for all 
	$ k \in \{ 1, \dots, L \} $, 
	$ i \in \{ 1, \ldots, \ell_k \} $,
	$ j \in \{ 1, \ldots, \ell_{ k - 1 } \} $ 
	that
	\begin{equation}
	\label{wb}
	\fw^{ k, \theta }_{ i, j }
	= 
	\theta_{ ( i - 1 ) \ell_{ k - 1 } + j 
		+ 
		\sum_{ h = 1 }^{ k - 1 } \ell_h ( \ell_{ h - 1 } + 1 ) }
	\qqandqq
	\fb^{ k, \theta }_i 
	=
	\theta_{ \ell_k \ell_{ k - 1 } + i 
		+ 
		\sum_{ h = 1 }^{ k - 1 } \ell_h ( \ell_{ h - 1 } + 1 ) } 
	,
	\end{equation}
	for every 
	$ k \in \N $, 
	$ \theta \in \R^{ \fd } $ 
	let 
	$
	\cA_k^{ \theta }
	=
	( 
	\cA_{ k, 1 }^{ \theta }, \ldots, \cA_{ k, \ell_k }^{ \theta } 
	)
	\colon \R^{ \ell_{ k - 1 } } \to \R^{ \ell_k } 
	$
	satisfy 
	for all 
	$ x \in \R^{ \ell_{ k - 1 } } $
	that 
	\begin{equation}
	\label{eq:def_Ak_transformation_deep_ANNs_normalized}
	\cA_k^{ \theta }( x ) 
	= 
	\fb^{ k, \theta } + \fw^{ k, \theta } x 
	,
	\end{equation}
	let 
	$ \Rect_r \colon \R \to \R $, 
	$ r \in \N \cup \{ \infty \} $,
	satisfy for all 
	$ 
	x \in \R
	$
	that 
	$
	(
	\cup_{ r \in \N }
	\{ \Rect_r \}
	)
	\subseteq C^1( \R, \R )
	$,
	$
	\sup_{ r \in \N }
	\sup_{ y \in [ - |x|, |x| ] }
	| ( \Rect_r )'( y ) |
	< \infty
	$,
	$
	\Rect_{ \infty }( x ) = \max\{ x, 0 \}
	$,
	and 
	\begin{equation}
	\label{lim_R}
	\textstyle
	\limsup_{ R \to \infty }
	\bigl[
	\sum_{ r = R }^{ \infty }
	\mathbbm{1}_{ (0,\infty) 
	}\bigl(
	| 
	\Rect_r( x ) - \Rect_{ \infty }( x ) 
	|
	+
	|
	( \Rect_r )'( x ) 
	-
	\mathbbm{1}_{ (0,\infty) }( x )
	|
	\bigr)
	\bigr]
	= 0
	,
	\end{equation}
	for every 
	$ r \in \N \cap \{ \infty \} $, 
	$ k \in \N $
	let 
	$ \mathfrak{M}_{ r }^{ k } \colon \R^{ \ell_k } \to \R^{ \ell_k } $ 
	satisfy for all 
	$ x = ( x_1, \ldots, x_{ \ell_k } ) \in \R^{ \ell_k } $ 
	that 
	\begin{equation}
	\textstyle
	\mathfrak{M}_{ r }^{ k }( x ) 
	= ( \Rect_r( x_1 ), \ldots, \Rect_r( x_{ \ell_k } ) )
	,
	\end{equation}
	let $ \mu \colon \cB( [ a, b ]^{ \ell_0 } ) \to [0, \infty] $ 
	be a finite measure, 
	for every 
	$ \theta \in \R^{ \fd } $, 
	$ r \in \N \cap \{ \infty \} $
	let 
	$
	\mathcal{N}^{ k, \theta }_{ r }
	= 
	( 
	\mathcal{N}^{ k, \theta }_{ r, 1 }, 
	\dots, 
	\mathcal{N}^{ k, \theta }_{ r, \ell_k } 
	)
	\colon \R^{ \ell_0 } \to \R^{ \ell_k } 
	$, 
	$ k \in \N $, 
	satisfy 
	for all 
	$ k \in \N $, 
	$ i \in \{ 1, \ldots, \ell_k \} $, 
	$ x \in \R^{ \ell_0 } $
	that
	$
	\mathcal{N}^{ 1, \theta }_r = \cA^{ \theta }_1( x )
	$
	and 
	\begin{equation}
	\label{eq:def_NN_realization_normalized}
	\textstyle 
	\mathcal{N}^{ k + 1, \theta }_{ r }( x )
	= 
	\cA_{ k + 1 }^{ \theta }\bigl(
	( 
	\mathfrak{M}_{ r}^{ k } \circ 
	\mathcal{N}^{ k, \theta }_{ r } 
	)( x )
	-
	\mathbbm{1}_{ 
		\{ L \}
	}( k + 1 )
	\int_{ [a,b]^{ \ell_0 } }
	( 
	\mathfrak{M}_{ r }^{ k } \circ 
	\mathcal{N}^{ k, \theta }_{ r } 
	)( y )
	\,
	\mu( \d y )
	\bigr)
	,
	\end{equation}
	let
	$ 
	f = 
	( f_1, \ldots, f_{ \ell_L } ) 
	\colon [a,b]^{ \ell_0 } \to \R^{ \ell_L } 
	$
	be measurable, 
	for every 
	$ r \in \N \cup \{ \infty \} $ 
	let $ \cL_{ r } \colon \R^{ \fd } \to \R $ 
	satisfy 
	for all $ \theta \in \R^{ \fd } $
	that
	\begin{equation}
	\label{eq:def_risk_function_normalized}
	\textstyle 
	\cL_{ r }( \theta ) 
	= 
	\int_{ [a,b]^{ \ell_0 } } 
	\norm{ \mathcal{N}_r^{ L, \theta }( x ) - f(x) }^2 \, \mu( \d x )
	,
	\end{equation}
	let 
	$
	\cG = ( \cG_1, \ldots, \cG_{ \fd } ) \colon \R^{ \fd } \to \R^{ \fd } 
	$ 
	satisfy
	for all 
	$
	\theta \in 
	\{
	\vartheta \in \R^{ \fd } \colon 
	( ( \nabla\cL_{ r } )( \vartheta ) )_{ r \in \N }
	\text{ is convergent} 
	\}
	$ 
	that
	$
	\cG( \theta ) = \lim_{ r \to \infty }( \nabla\cL_{ r } )( \theta ) 
	$, 
	for every 
	$ k \in \N $,
	$ \theta \in \R^{ \fd } $, 
	$ i \in \{ 1, \dots, \ell_k \} $
	let 
	$ 
	V^{ k, \theta }_i
	=
	(
	V^{ k, \theta }_{ i, 1 }
	,
	\dots ,
	V^{ k, \theta }_{ i, \ell_{ k - 1 } + 1 }
	)
	\in \R^{ 
		\ell_{ k - 1 } + 1
	}
	$
	satisfy 
	\begin{equation}
	V^{ k, \theta }_i
	=
	(
	\fw^{ k, \theta }_{ i, 1 }
	,
	\dots, 
	\fw^{ k, \theta }_{ i, \ell_k }
	,
	\fb^{ k, \theta }_i
	)
	,
	\end{equation}
	for every $ k \in \N $, $ i \in \{ 1, \dots, \ell_k \} $ 
	let 
	$
	\psi^k_i \colon \R^{ \fd } \to \R
	$
	satisfy for all $ \theta \in \R^{ \fd } $ that
	$
	\psi^k_i( \theta ) = \| V^{ k, \theta }_i \|^2
	$, 
	let 
	$
	\Lambda \subseteq \N^2 
	$
	satisfy 
	$
	\Lambda = 
	\cup_{ k = 1 }^{ L - 1 }
	(
	\{ k \} \times 
	\{ 1, 2, \dots, \ell_k \}
	)
	$, 
	let 
	$ 
	\rho \colon 
	( \cup_{ n \in \N } \R^n ) 
	\to 
	( \cup_{ n \in \N } \R^n ) 
	$
	satisfy for all 
	$ n \in \N $, $ x \in \R^n $ that 
	\begin{equation}
	\label{eq:setting:def:rho}
	\rho( x ) 
	=
	\bigl[ 
	\| x \| + \mathbbm{1}_{ \{ 0 \} }( \| x \| ) 
	\bigr]^{ - 1 } x
	,
	\end{equation}
	let 
	$ G \colon \R^{ \fd } \to \R^{ \fd } $
	satisfy for all 
	$ \theta \in \R^{ \fd } $ 
	that
	\begin{equation}
	\begin{split}
	\textstyle 
	G( \theta )
	=
	\cG( \theta )
	-
	\sum_{ (k,i) \in \Lambda }
	\bigl\langle
	\rho\bigl( 
	( \nabla \psi^k_i )( \theta ) 
	\bigr) 
	,
	\cG( \theta )
	\bigr\rangle
	\,
	\rho\bigl( 
	( \nabla \psi^k_i )( \theta ) 
	\bigr) 
	,
	\end{split}
	\end{equation} 
	let 
	$
	\phi \colon \R^{ \fd } \to \R^{ \fd }
	$
	satisfy 
	for all 
	$ k \in \{ 1, \dots, L \} $,  
	$ \theta \in \R^{ \fd } $, 
	$ i \in \{ 1, \dots, \ell_k \} $
	that 
	\begin{equation}
	V^{ k, \phi( \theta ) }_i
	= 
	\begin{cases}
	\rho( V^{ k, \theta }_i )
	&
	\colon 
	k < L
	\\[1ex]
	V^{ k, \theta }_i 
	&
	\colon 
	k = L
	,
	\end{cases}
	\end{equation}
let 
$
\Psi_k \colon \R^{ \fd } \to \R^{ \fd }
$, 
$ k \in \N_0 $, 
satisfy for all 
$ k, K \in \N $, 
$ i \in \{ 1, \dots, \ell_K \} $, 
$ \theta \in \R^{ \fd } $ 
that 
$
\Psi_0( \theta ) = \theta 
$
and 
\begin{equation}
\label{eq:setting:def:psik}
V^{ K, \Psi_k( \theta ) }_i
=
\begin{cases}
\rho(
V^{ K, \Psi_{ k - 1 }( \theta ) }_i
)
&
\colon
K = k 
\\[1ex]
\operatorname{diag}\bigl( 
\| V^{ k, \Psi_{ k - 1 }( \theta ) }_1 \|
,
\dots 
,
\| V^{ k, \Psi_{ k - 1 }( \theta ) }_{ \ell_{ k } } \|
, 
1
\bigr)
\,    
V^{ K, \Psi_{ k - 1 }( \theta ) }_i
&
\colon
K = k + 1
\\[1ex]
V^{ K, \Psi_{ k - 1 }( \theta ) }_i
&
\colon
K \notin \{ k, k + 1 \}
,
\end{cases}
\end{equation}
	let $ ( \Omega, \cF, \P ) $ be a probability space, 
	let $ \xi \colon \Omega \to \R^{ \fd } $ be a random variable, 
	let $ \Theta \colon [0,\infty) \times \Omega \to \R^{ \fd } $ 
	satisfy for all $ t \in [0,\infty) $, $\omega \in \Omega$ that 
	$
	\int_0^t
	\|  
	G( \Theta_s( \omega ) )
	\|
	\,
	\d s
	< \infty 
	$
	and 
	\begin{equation}
	\label{eq:setting:def:flow}
	\textstyle 
	\Theta_t( \omega ) = \Psi_{ L - 1 } ( \xi( \omega ) ) - \int_0^t G( \Theta_s( \omega ) ) \, \d s
	,
	\end{equation}
	let $ ( \gamma_n )_{ n \in \N_0 } \subseteq \R $,   
	and 
	let 
	$
	\varTheta_n \colon \Omega \to \R^{ \fd } 
	$, 
	$ n \in \N_0 $, 
	satisfy 
	for all $ n \in \N_0 $, $\omega \in \Omega$ that 
	$
	\varTheta_0( \omega ) = 
	\Psi _ { L - 1 } ( 
	\xi( \omega )
	)
	$
	and 
	\begin{equation}
	\label{eq:gradient_descent_def2}
	\varTheta_{ n + 1 }( \omega )
	= 
	\phi\bigl(
	\varTheta_n( \omega )
	-
	\gamma_n
	G( \varTheta_n( \omega ) )
	\bigr)
	.
	\end{equation}
\end{setting}

Next, in \cref{prop:rescaling:properties} we verify some basic properties of the rescaling operation $\Psi_{L - 1 } \colon \R^\fd \to \R^\fd$.
 In particular, we show for every parameter vector $\theta \in \cu{ \vartheta \in \R^\fd \colon \min_{ (k,i) \in \Lambda } \norm{V_i^{ k , \vartheta } } > 0 }$ that the rescaled vector $\Psi_{L - 1 } \in \R^\fd$ is an element of a suitable $C^\infty$-submanifold of the parameter space $\R^\fd$; see \cref{prop:properties:item1,prop:properties:item3}. In addition, we demonstrate that the rescaling map $\Psi_{L - 1 }$ does not change the output of the considered ANN with ReLU activation; see \cref{prop:properties:item2}.

\begin{prop} [Properties of ANNs with normalized parameter vectors]
	\label{prop:rescaling:properties}
	Assume \cref{setting:dnn_normalized}. 
	Then 
	\begin{enumerate}[label=(\roman*)]
		\item 
		\label{prop:properties:item1}
		it holds that 
		\begin{equation}
		\{ 
		\theta \in \R^{ \fd } 
		\colon 
		(
		\forall \, (k,i) \in \Lambda 
		\colon
		\| V^{ k, \theta }_i \| = 1
		)
		\}
		\end{equation}
		is a 
		$ ( \fd - ( \sum_{ k = 1 }^{ L - 1 } \ell_k ) ) $-dimensional 
		$ C^{ \infty } $-submanifold 
		of the $ \R^{ \fd } $, 
		\item 
		\label{prop:properties:item2}
		it holds for all $ \theta \in \R^{ \fd } $ that
		\begin{equation}
		\cN^{ L, \theta }_{ \infty }
		=
		\cN^{ L, \Psi _ { L - 1 } ( \theta ) }_{ \infty }
		,
		\end{equation}
		and
		\item
		\label{prop:properties:item3}
		 it holds for all $\theta \in \cu{ \vartheta \in \R^\fd \colon \min_{ (k,i) \in \Lambda } \norm{V_i^{ k , \vartheta } } > 0 }$,
		 $ (k , i ) \in \Lambda$ that $ \norm{ V_i^{ k , \Psi_{L- 1 } (  \theta  ) } } = 1 $.
	\end{enumerate}
\end{prop} 

\begin{cproof}{prop:rescaling:properties}
First, to \prove[s] \cref{prop:properties:item1} let $\cV \colon \R^\fd \to \R^{ \# \Lambda }$
	satisfy for all $\theta \in \R^\fd$ that $\cV ( \theta ) = (\psi_i^k ( \theta ) )_{ (k , i ) \in \Lambda }$.
	\Nobs that $\cV \in C ^\infty ( \R^\fd , \R^{ \# \Lambda })$.
	In addition, \nobs that for all $\theta \in \cu{\vartheta \in \R^\fd \colon ( \forall \, (k , i ) \in \Lambda \colon \psi_i^k ( \vartheta ) \not= 0 ) }$ it holds that $\operatorname{rank} ( \cV' ( \theta ) ) = \# \Lambda = \sum_{k=1}^{L - 1 } \ell_k$.
	Combining this with the preimage theorem (cf., e.g., Tu~\cite[Theorem 9.9]{Tu2011})
	\proves that
	\begin{equation}
	\cu[\big]{ \theta \in \R^{ \fd } 
		\colon 
		( \forall \, (k,i) \in \Lambda 
		\colon
		\norm{ V^{ k, \theta }_i } = 1
		) } 
	= \cu{\theta \in \R^\fd \colon \cV ( \theta ) = (1, 1, \ldots, 1 ) }
	\end{equation}
	is a $ ( \fd - ( \sum_{ k = 1 }^{ L - 1 } \ell_k ) ) $-dimensional 
	$ C^{ \infty } $-submanifold 
	of the $ \R^{ \fd } $. This \proves[e] \cref{prop:properties:item1}.
	
	To \prove[p] \cref{prop:properties:item2} let $K \in \cu{1, 2, \ldots, L - 1 }$, $\theta \in \R^\fd$, $y \in \R^{ \ell_0}$ be fixed
	and for every $\vartheta \in \R^\fd$ let $\cN_\infty ^{ 0 , \vartheta } \colon \R^{ \ell_0 } \to \R^{ \ell_0}$ satisfy for all $x \in \R^{ \ell_0}$ that $\cN_\infty ^{ 0 , \vartheta } ( x ) = x $.
	\Nobs that the fact that for all $k \in \cu{1, 2, \ldots, K - 1 }$,
	$i \in \cu{1, 2, \ldots, \ell_k}$ it holds that $V_i^{ k , \Psi_K ( \theta ) } = V_i^{ k , \Psi_{ K - 1 } ( \theta ) }$ \proves that $\cN_\infty ^{ K - 1 , \Psi_K ( \theta ) } ( y ) = \cN_\infty ^{ K - 1 , \Psi_{K - 1 } ( \theta ) } ( y ) $.
	Moreover,
	\nobs that \cref{eq:setting:def:psik} \proves for all $i \in \cu{1, 2 , \ldots, \ell_K}$, $j \in \cu{1, 2, \ldots, \ell_{K - 1 } } $ that
	\begin{equation}
	\begin{split}
	\fw_{i,j} ^{ K , \Psi_K ( \theta ) } 
	& =  \fw_{i,j}^{ K , \Psi_{K - 1 } ( \theta ) } \rbr[\big]{\norm{ V_i^{ K , \Psi_{K - 1 } ( \theta )  } } + \indicator{ \cu{0 } } ( \norm{ V_i^{ K , \Psi_{K - 1 } ( \theta )  } } ) }^{ - 1 } , \\
	\fb_{i} ^{K , \Psi_K ( \theta ) }
	& = \fb_{i}^{ K , \Psi_{K - 1 } ( \theta ) } \rbr[\big]{\norm{ V_i^{ K , \Psi_{K - 1 } ( \theta )  } } + \indicator{ \cu{0 } } ( \norm{ V_i^{ K , \Psi_{K - 1 } ( \theta )  } } ) } ^{ - 1 } .
	\end{split}
	\end{equation}
	Therefore, we get for all $i \in \cu{1, 2, \ldots, \ell_{K } }$ that
	\begin{equation}
	\begin{split}
	&\cN_{\infty , i } ^{ K , \Psi_K ( \theta ) } ( y )
	= \fb_{i} ^{ K , \Psi_K ( \theta ) }  + \ssuml_{j=1}^{ \ell_{K - 1 } } \fw_{i,j} ^{ K , \Psi_K ( \theta ) } 
	\Rect_{ \infty } \rbr[\big]{ \cN_{\infty , j } ^{ \bfk - 1 , \Psi_{K } ( \theta ) } ( y ) }
	\\
	&= 
	\rbr[\bigg]{ \fb_{i} ^{ K , \Psi_{K - 1 } ( \theta ) } + \ssuml_{j=1}^{ \ell_{K - 1 } } \fw_{i,j} ^{ K, \Psi_{ K - 1 } ( \theta ) } 
		\Rect_{ \infty } \rbr[\big]{ 	\cN_{\infty , j } ^{ K - 1 , \Psi_{K - 1 } ( \theta ) } ( y ) } } \\
	& \qquad \times \rbr[\big]{ \norm{ V_i^{ K , \Psi_{K - 1 } ( \theta )  } } + \indicator{ \cu{0 } } ( \norm{ V_i^{ K , \Psi_{K - 1 } ( \theta )  } } ) } ^{ - 1} \\
	&= \rbr[\big] { \norm{ V_i^{ K , \Psi_{K - 1 } ( \theta )  } } + \indicator{ \cu{0 } } ( \norm{ V_i^{ K , \Psi_{K - 1 } ( \theta )  } } ) } ^{-1}
	 \cN_{\infty , i } ^{ K , \Psi_{ K - 1 } ( \theta ) } ( y ) .
	\end{split}
	\end{equation}
	This and the fact that $\forall \, u \in \R , \, \eta \in [0 , \infty ) \colon \Rect_{ \infty } ( \eta u ) = \eta \Rect_{ \infty } ( u ) $
	\prove for all $i \in \cu{1, 2, \ldots, \ell_{K } }$ that
	\begin{equation}
	\begin{split}
	\Rect_{ \infty } ( \cN_{\infty , i } ^{ K , \Psi_K ( \theta ) } ( y ) ) 
	 =  \rbr[\big]{ \norm{ V_i^{ K , \Psi_{K - 1 } ( \theta )  } } + \indicator{ \cu{0 } } ( \norm{ V_i^{ K , \Psi_{K - 1 } ( \theta )  } } ) } ^{-1}
	  \Rect_{ \infty } ( \cN_{\infty , i } ^{ K , \Psi_{ K - 1 } ( \theta ) } ( y ) ).
	\end{split} 
	\end{equation}
	In addition,
	\nobs that \cref{eq:setting:def:psik} \proves for all 	
	$i \in \cu{1, 2 , \ldots, \ell_{K + 1 } }$, $j \in \cu{1, 2, \ldots, \ell_{K} } $ that
	\begin{equation}
	\fw_{i,j} ^{ K + 1 , \Psi_{K } ( \theta ) } = \fw_{i,j} ^{ K + 1 , \Psi_{K - 1 } ( \theta ) } \norm{ V_j^{ K , \Psi_{K - 1 } ( \theta ) } }
	\qandq 
	\fb_i^{ K + 1 , \Psi_K ( \theta ) } = \fb_i^{ K + 1 , \Psi_{K - 1 } ( \theta ) }.
	\end{equation}
	Combining this with the fact that for all
	$j \in \cu{1, 2, \ldots, \ell_K }$ with $\norm{ V_j^{ K , \Psi_{K - 1 } ( \theta ) } } = 0$
	it holds that $\cN_{\infty , j } ^{ K , \Psi_{ K - 1 } ( \theta ) } ( y ) = \cN_{\infty , j } ^{ K , \Psi_K ( \theta ) } ( y ) = 0$
	\proves that for all 
	$i \in \cu{1, 2 , \ldots, \ell_{K + 1 } }$
	we have that
	\begin{equation}
	\begin{split}
	&\fb_i^{ K + 1 , \Psi_K ( \theta ) } + \ssuml_{j=1}^{ \ell_K } \fw_{i,j} ^{ K + 1 , \Psi_{K } ( \theta ) } \Rect_{ \infty } ( \cN_{\infty , j } ^{ K , \Psi_K ( \theta ) } ( y ) ) \\
	&= \fb_i^{ K + 1 , \Psi_{K - 1 } ( \theta ) } +
	\ssuml_{j=1}^{ \ell_K } \tfrac{ \norm{ V_j^{ K , \Psi_{K - 1 } ( \theta ) } } } { \norm{ V_j^{ K , \Psi_{K - 1 } ( \theta )  } } + \indicator{ \cu{0 } } ( \norm{ V_j ^{ K , \Psi_{K - 1 } ( \theta )  } } ) } 
	\fw_{i,j} ^{ K + 1 , \Psi_{K - 1 } ( \theta ) } \Rect_{ \infty } ( \cN_{\infty , j } ^{ K , \Psi_{ K - 1 } ( \theta ) } ( y ) ) \\
	&= \fb_i^{ K + 1 , \Psi_{K - 1 } ( \theta ) } +
	\ssuml_{j=1}^{ \ell_K } 
	\fw_{i,j} ^{ K + 1 , \Psi_{K - 1 } ( \theta ) } \Rect_{ \infty } ( \cN_{\infty , j } ^{ K , \Psi_{ K - 1 } ( \theta ) } ( y ) ) .
	\end{split}
	\end{equation}
	This and \cref{eq:def_NN_realization_normalized} \prove[dpsin] that $\cN_\infty ^{ K + 1  , \Psi_K ( \theta ) } ( y ) = \cN_\infty ^{ K + 1 , \Psi_{K - 1 } ( \theta ) } ( y ) $.
	Hence, we obtain that $\cN_\infty ^{ L  , \Psi_K ( \theta ) } = \cN_\infty ^{ L , \Psi_{K - 1 } ( \theta ) } $.
	Induction therefore \proves[e] \cref{prop:properties:item2}.
	
	Next \nobs that \cref{eq:setting:def:psik},
	\cref{eq:setting:def:rho},
	and induction \prove for all
	$\vartheta \in \R^\fd $, $j \in \cu{1, 2, \ldots, L - 1 }$, $(k , i ) \in \Lambda $ with $ \norm{V_i^{ k , \vartheta' } } > 0$ and $k \le j$ that 
	$\norm{ V_i^{ k , \Psi_j ( \vartheta ) } } = 1$.
	This \proves[e] \cref{prop:properties:item3}.
\end{cproof} 

In the following result, \cref{prop:modified:gf:properties}, we establish some invariance properties of the considered modified GF and GD processes in \cref{setting:dnn_normalized}. 
In particular, we show for every $\omega \in \Omega$ for which the initial value $\xi ( \omega ) \in \R^\fd$ is non-degenerate in a suitable sense that the corresponding GF trajectory $ (\Theta_t ( \omega ) ) _ { t \in [ 0 , \infty ) }$ stays on the considered $C^\infty$-submanifold of the parameter space (see \cref{prop:properties:item4})
and has non-increasing risk value $\cL _\infty ( \Theta_t ( \omega ) ) $, $t \in [0 , \infty )$ (see \cref{prop:properties:item5}).
For the proof we employ \cref{lem:modified_gradient} and \cref{cor:descent}.

\begin{prop} [Properties of modified GF and GD processes]
	\label{prop:modified:gf:properties}
	Assume \cref{setting:dnn_normalized}. 
	Then 
	\begin{enumerate}[label=(\roman*)]
		\item 
		\label{prop:properties:item4}
		it holds 
		for all 
		$ (k,i) \in \Lambda $, 
		$ t \in [0,\infty) $, $\omega \in \Omega$
		with 
		$ 
		\min_{
			(k,i) \in \Lambda
		}
		\| V^{ k, \xi( \omega ) }_i \|
		> 0
		$
		that 
		\begin{equation}
		\psi^k_i( \Theta_t( \omega ) ) 
		= 
		1
		,
		\end{equation}
		\item 
		\label{prop:properties:item5}
		it holds for all 
		$ s \in [0,\infty) $,
		$t \in [s , \infty )$,
		 $\omega \in \Omega$
		with 
		$ 
		\min_{
			(k,i) \in \Lambda
		}
		\| V^{ k, \xi( \omega ) }_i \|
		> 0
		$
		that $\cL_\infty ( \Theta_t ( \omega ) ) \le \cL_\infty ( \Theta_s ( \omega ) )$,
		and 
		\item 
		\label{prop:properties:item6}
		it holds 
		for all 
		$ (k,i) \in \Lambda $, 
		$ n \in \N_0 $, $\omega \in \Omega$ 
		with 
		$
		\min_{ (k,i) \in \Lambda }
		\|
		V^{ k, 
			\varTheta_n( \omega ) - \gamma_n G( \varTheta_n( \omega ) )
		}_i
		\|
		> 0
		$
		that 
		\begin{equation}
		\psi^k_i( \varTheta_{ n + 1 }( \omega ) )
		= 
		1
		.
		\end{equation}
	\end{enumerate}
\end{prop}

\begin{cproof}{prop:modified:gf:properties}
First, to prove \cref{prop:properties:item4} let $\omega \in \Omega$ satisfy 	
$ \min_{ (k,i) \in \Lambda } 
\norm{ V^{ k, \xi( \omega ) }_i }
> 0 $
and denote
\begin{equation}
\label{prop:properties:eq:deftau}
\tau = \inf \rbr[\big]{ \cu[\big]{ t \in [0 , \infty ) \colon \min\nolimits_{ (k,i) \in \Lambda } \norm{ V_i^{ k , \Theta_t ( \omega ) } } = 0} \cup \cu{ \infty } } \in [0 , \infty ].
\end{equation}
 \Nobs that \cref{prop:properties:item3} in \cref{prop:rescaling:properties} and \cref{eq:setting:def:flow} \prove for all $(k , i ) \in \Lambda$ that $ \norm{ V_i^{ k , \Theta_0 ( \omega ) } } = \norm{  V_i^{ k , \Psi_{L - 1 } ( \xi ( \omega ) ) } } = 1 $. 
 Hence, we obtain for all $(k , i) \in \Lambda$ that $\psi_i^k ( \Theta_0 ( \omega ) ) = 1$.
 Furthermore, the fact that $\Theta$ is continuous \proves that $\tau > 0$.
 In addition, \nobs that for all $t \in [0 , \tau )$ we have that
 \begin{equation}
 \begin{split}
 G ( \Theta_t ( \omega ) ) 
 &= \cG ( \Theta_t ( \omega )  )
 - \ssum_{ ( k , i ) \in \Lambda } \tfrac{ \spro{ \nabla \psi_i^k ( \Theta_t ( \omega ) ) , \cG ( \Theta_t ( \omega ) ) } } { \norm{ \nabla \psi_i^k ( \Theta_t ( \omega ) ) } ^2 } \nabla \psi_i^k ( \Theta_t ( \omega ) ).
 \end{split}
 \end{equation}
 Combining this with
 \cref{eq:setting:def:flow}, the fact that for all $(k , i ) \in \Lambda$ it holds that $\psi_i^k \in C^1 ( \R^\fd , \R )$,
 the fact that for all $t \in [0 , \tau )$ it holds that
 $\inf_{s \in [0 , t ] } \min_{ ( k , i ) \in \Lambda } \norm{ \nabla \psi_i^k ( \Theta_s ( \omega ) ) } > 0$,
 and \cref{lem:modified_gradient}
 (applied with $K \with \# \Lambda $, $(\psi_i)_{i \in \cu{1, 2, \ldots, K } } \with ( \psi_i^k)_{ ( k , i ) \in \Lambda } $ in the notation of \cref{lem:modified_gradient})
 \proves for all $t \in [0 , \tau )$
 that
 $\psi^k_i( \Theta_ t ( \omega ) ) = \norm{ V_i^{ k , \Theta_t ( \omega ) } } ^2 = 1$.
 This, \cref{prop:properties:eq:deftau},
  and the fact that $\Theta \in C( [ 0 , \infty ) , \R^\fd )$ \prove that $\tau = \infty $,
 which \proves[e] \cref{prop:properties:item4}.
 
 Next \nobs that Hutzenthaler et al.~\cite[Theorem 2.9]{HutzenthalerJentzenPohlRiekertScarpa2021} \proves for all $\theta \in \R^\fd$ that $\bigcup_{r \in \N } \cu{ \cL_r } \subseteq C^1 ( \R^\fd , \R )$, $\lim_{r \to \infty } \cL_r ( \theta ) = \cL_\infty ( \theta )$,
  and $\lim_{r \to \infty } \nabla \cL_r ( \theta ) = \cG ( \theta )$.
 Furthermore, \cite[Lemma 3.6]{HutzenthalerJentzenPohlRiekertScarpa2021} \proves for all compact $K \subseteq \R^\fd$ that $\sup_{\theta \in K } \sup_{r \in \N } \norm{ \nabla \cL_r ( \theta ) } < \infty $.
 Combining this with \cref{cor:descent} \proves[e] \cref{prop:properties:item5}.
 
Finally, to \prove[p] \cref{prop:properties:item6} let
 $n \in \N_0$,
  $\omega \in \Omega$ satisfy
  $
  \min_{ (k,i) \in \Lambda }
  \norm{ V^{ k,  \varTheta_n( \omega ) - \gamma_n G( \varTheta_n( \omega ) )  }_i}
   > 0$.
   \Nobs that \cref{eq:setting:def:rho}
   \proves for all
   $(k , i ) \in \Lambda$
   that
   \begin{equation}
   \begin{split}
   \psi_i^k ( \varTheta _{n+1} ( \omega ) ) 
   &= \norm{ V_i^{ k , \varTheta _{n+1} ( \omega ) } } ^2 
   = \norm[\big]{ \rho \rbr[\big]{ V_i^{ k , \varTheta_n ( \omega ) - \gamma_n G ( \varTheta_n ( \omega ) ) } } } ^2 = 1 .
   \end{split}
   \end{equation}
   This \proves[e] \cref{prop:properties:item6}.
\end{cproof}

\section{Global boundedness of normalized gradient flows in the training of shallow ReLU ANNs with one hidden neuron}
\label{sec:one_neuron_analysis}

In this section we prove that the modified gradient flow considered in \cref{sec:normalized_gradient_flow} is uniformly bounded in the case of shallow ANNs with one-dimensional input, one neuron on the hidden layer,
 one-dimensional output,
and uniformly distributed input data; see \cref{theo:bounded:1neuron} below. For convenience we first introduce the simplified notation we will employ throughout this section.

\subsection{Notation} \label{subsec:1n:notation}
Let $\mu \colon \cB ( \R ) \to [0 , \infty ]$ be the Lebesgue measure on $\R$.
Let $f \in C (  [0,1] , \R )$ be the target function and
let $\overline{f} = \int_0^1 f(x) \, \d x$.
Let $m \colon \R^3 \to \R$
satisfy for every $\theta = ( \theta_1, \theta_2 , \theta_3 ) \in \R^3$
that $m ( \theta ) = \int_0^1  \max \cu{\theta_1 s + \theta_2 , 0 }  \, \d s$.
We consider the risk function\footnote{For simplicity we assume that the outer bias has the fixed value $\overline{f}$. So the risk function only depends on three parameters.}
$\cL \colon \R^3 \to \R$ which satisfies for all $\theta \in \R^3$ that
\begin{equation}
\cL ( \theta ) = \int_0^1 \rbr*{ \theta_3 ( \max \cu{\theta_1 s + \theta_2 , 0 } - m ( \theta ) ) + \overline{f} - f ( s ) } ^2 \ \d s.
\end{equation}
Let $g \colon \R^3 \to \R$ satisfy for all $\theta \in \R^3$ that $g ( \theta ) = \abs{\theta_1 }^2 + \abs{\theta_2 }^2$
and consider the two-dimensional $C^\infty$-manifold $\cM = g^{-1} ( 1 ) \subseteq \R^3$.
We want the gradient flow to stay on this manifold.

\Nobs that Ibragimov et al.~\cite[Corollary 2.3]{IbragimovLocalMinima} ensures for all $\theta \in \R^3$ with $\abs{\theta_1} + \abs{\theta_2 } > 0$ that $m$ is differentiable at $\theta$.
Using \cite[Corollary 2.3]{IbragimovLocalMinima} again ensures for all $\theta \in \R^3$ with $\abs{\theta_1} + \abs{\theta_2 } > 0$
that all partial derivatives of $\cL$ at $\theta$ exist. It is also not hard to see that these derivatives are continuous around $\theta$. 
Consider a modified gradient $\fG \colon \R^3 \to \R^3$ which is locally bounded and measurable and satisfies for all $\theta \in \R^3$ with $\abs{\theta_1} + \abs{\theta_2 } > 0$ that
\begin{equation}
\fG ( \theta ) = \nabla \cL ( \theta ) - \norm{ \nabla g ( \theta ) } ^{ - 2 }  \spro{\nabla \cL ( \theta ) , \nabla g ( \theta ) } \nabla g ( \theta ) .
\end{equation}
From \cite[Proposition 2.11]{JentzenRiekertFlow} we know for all
 $\theta \in \R^3$ with $\abs{\theta_1} + \abs{\theta_2 } > 0$ that $\nabla \cL ( \theta )$ agrees with the generalized gradient introduced in \cref{setting:dnn_normalized}.
For every $\theta \in \R^3$ let $I^\theta = \cu{s \in [0 , 1 ] \colon \theta_1 s + \theta_2 > 0 }$
and let \begin{equation}
q^\theta = \begin{cases}
- \tfrac{\theta_2 }{ \theta_1 } & \colon \theta_1 \not= 0 \\
\infty & \colon \text{else} .
\end{cases}
\end{equation}
In the following we consider a gradient flow (GF) trajectory $\Theta = ( \Theta_1 , \Theta_2 , \Theta_3 ) \colon [0, \infty ) \to \R^3$ which satisfies $\Theta ( 0 ) \in \cM $ and
$\forall \, t \in [0 , \infty ) \colon \Theta ( t ) = \Theta ( 0 ) - \int_0^t \fG  ( \Theta ( u ) ) \, \d u$.

\subsection{Basic properties of the gradient and the GF trajectory}

\begin{lemma}
	\label{lem:gradient:basic}
	Consider the notation in \cref{subsec:1n:notation}. Then
	\begin{enumerate} [label = (\roman*)]
		\item
		\label{lem:gradient:basic:item1}
		it holds for all $\theta \in \cM$ that
		\begin{equation}
		\label{eq:gen:gradient:explicit}
		\begin{split}
		\fG_1 ( \theta )
		& = 2 \theta_3  \int_0^1 \rbr*{ \theta_3 ( \max \cu{\theta_1 s + \theta_2 , 0 } - m ( \theta ) ) + \overline{f} - f ( s ) } ( \theta_2^2 s - \theta_1 \theta_2 ) \indicator{I^\theta} ( s ) \, \d s, \\
		\fG_2 ( \theta )
		& = 2 \theta_3 \int_0^1 \rbr*{ \theta_3 ( \max \cu{\theta_1 s + \theta_2 , 0 } - m ( \theta ) ) + \overline{f} - f ( s ) } ( \theta_1 ^2 - \theta_1 \theta_2 s ) \indicator{I^\theta} ( s )  \, \d s, \\
		\fG_3 ( \theta)
		&= 2 \int_0^1 \rbr*{ \theta_3 ( \max \cu{\theta_1 s + \theta_2 , 0 } - m ( \theta ) ) + \overline{f} - f ( s ) }\\
		& \qquad \times \rbr*{ \max \cu{\theta_1 s + \theta_2 , 0 } - m ( \theta ) } \, \d s 
		\end{split}
		\end{equation}
		and 
		\item
		\label{lem:gradient:basic:item2}
		it holds that $\fG  | _\cM \colon \cM  \to \R^3$ is locally Lipschitz continuous.
	\end{enumerate}
\end{lemma}

\begin{cproof}{lem:gradient:basic}
	First, \cite[Corollary 2.3]{IbragimovLocalMinima} \proves for all $\theta \in \R^3$ with $\abs{\theta_1} + \abs{\theta_2 } > 0$ that $m$ is differentiable at $\theta$ and satisfies $\frac{\partial}{\partial \theta_1} m ( \theta ) = \int_{I^\theta  } s \, \d s$ and $\frac{\partial}{\partial \theta_2} m ( \theta ) = \int_{I^\theta } 1 \, \d s$.
	This, \cite[Corollary 2.3]{IbragimovLocalMinima},
	and the fact that $\forall \, \theta \in \R^3 \colon \int_0^1 \rbr{ \theta_3 ( \max \cu{\theta_1 s + \theta_2 , 0 } - m ( \theta ) ) + \overline{f} - f ( s ) } \, \d s  = 0$ \prove for all $\theta \in \R^3$ with $\abs{\theta_1} + \abs{\theta_2 } > 0$ that
	\begin{equation} \label{lem:gradient:basic:eq1}
	\begin{split}
	\tfrac{\partial}{\partial \theta_1} \cL ( \theta )
	& = 2 \theta_3 \int_0^1 \rbr*{ \theta_3 ( \max \cu{\theta_1 s + \theta_2 , 0 } - m ( \theta ) ) + \overline{f} - f ( s ) }
	\rbr*{   s \indicator{I^\theta} ( s ) - \tint_{I^\theta} r \, \d r } \, \d s \\
	&= 2 \theta_3 \int_0^1 \rbr*{ \theta_3 ( \max \cu{\theta_1 s + \theta_2 , 0 } - m ( \theta ) ) + \overline{f} - f ( s ) } s \indicator{I^\theta } ( s ) \, \d s , \\
	\tfrac{\partial}{\partial \theta_2} \cL ( \theta )
	& = 2 \theta_3 \int_0^1 \rbr*{ \theta_3 ( \max \cu{\theta_1 s + \theta_2 , 0 } - m ( \theta ) ) + \overline{f} - f ( s ) }
	\rbr*{  \indicator{I^\theta} ( s ) - \tint_{I^\theta} 1 \, \d r } \, \d s \\
	&= 2 \theta_3 \int_0^1 \rbr*{ \theta_3 ( \max \cu{\theta_1 s + \theta_2 , 0 } - m ( \theta ) ) + \overline{f} - f ( s ) } \indicator{I^\theta } ( s ) \, \d s ,  \\
	\tfrac{\partial}{\partial \theta_3} \cL ( \theta )
	& = 2 \int_0^1 \rbr*{ \theta_3 ( \max \cu{\theta_1 s + \theta_2 , 0 } - m ( \theta ) ) + \overline{f} - f ( s ) }
	\rbr*{ \max \cu{\theta_1 s + \theta_2 , 0 } - m ( \theta ) } \, \d s .
	\end{split}
	\end{equation}
	Combining this with the fact that $\forall \, \theta \in \R^3 \colon \nabla g ( \theta ) = ( 2 \theta_1 , 2 \theta_2 , 0 ) ^T$
	\proves for all
	$\theta \in \R^3$ with $\abs{\theta_1} + \abs{\theta_2 } > 0$ that
	\begin{equation}
	\begin{split}
	\spro{ \nabla \cL ( \theta ) , \nabla g ( \theta ) }
	&= 4 \theta_3 \int_0^1 \rbr*{ \theta_3 ( \max \cu{\theta_1 s + \theta_2 , 0 } - m ( \theta ) ) + \overline{f} - f ( s ) } \\
	& \quad \times
	\rbr*{  \theta_1 s \indicator{I^\theta} ( s )+ \theta_2 \indicator{I^\theta} ( s )  } \, \d s , \\
	&= 4 \theta_3 \int_0^1 \rbr*{ \theta_3 ( \max \cu{\theta_1 s + \theta_2 , 0 } - m ( \theta ) ) + \overline{f} - f ( s ) } 
	\\ & \quad \times
	\rbr*{\theta_1 s + \theta_2 } \indicator{I^\theta } ( s )  \, \d s .
	\end{split}
	\end{equation}
	In addition, \nobs that for all $\theta \in \cM$ we have that $\norm{\nabla g ( \theta ) } ^2 = (2 \theta_1 )^2 + (2 \theta_2 ) ^2 = 4$.
	Therefore, we obtain for all $\theta \in \cM$
	that
	\begin{equation}
	\begin{split}
	\fG_1 ( \theta ) &=
	2 \theta_3 \int_0^1 \rbr*{ \theta_3 ( \max \cu{\theta_1 s + \theta_2 , 0 } - m ( \theta ) ) + \overline{f} - f ( s ) } ( s - \theta_1 ( \theta_1 s + \theta_2 ) ) \indicator{I^\theta} ( s ) \, \d s \\
	&=  2 \theta_3 \int_0^1 \rbr*{ \theta_3 ( \max \cu{\theta_1 s + \theta_2 , 0 } - m ( \theta ) ) + \overline{f} - f ( s ) } ( \theta_2^2 s - \theta_1  \theta_2 ) \indicator{I^\theta} ( s ) \, \d s , \\
	\fG_2 ( \theta ) &=
	2 \theta_3  \int_0^1 \rbr*{ \theta_3 ( \max \cu{\theta_1 s + \theta_2 , 0 } - m ( \theta ) ) + \overline{f} - f ( s ) } (1 - \theta_2 ( \theta_1 s + \theta_2 ) ) \indicator{I^\theta} ( s ) \, \d s \\
	&=  2 \theta_3 \int_0^1 \rbr*{ \theta_3 ( \max \cu{\theta_1 s + \theta_2 , 0 } - m ( \theta ) ) + \overline{f} - f ( s ) } ( \theta_1^2 - \theta_1  \theta_2 s ) \indicator{I^\theta} ( s ) \, \d s ,
	\end{split}
	\end{equation}
	and $\fG_3 ( \theta ) = \tfrac{\partial}{\partial \theta_3} \cL ( \theta )$.
	Combining this with \cref{lem:gradient:basic:eq1} establishes \cref{lem:gradient:basic:item1}.
	\Cref{lem:gradient:basic:item2} follows from \cite[Corollary 2.6]{EberleJentzenRiekert2021}.
\end{cproof}

As a consequence of \cref{lem:gradient:basic:item2},
\cref{lem:modified_gradient},
and \cref{cor:descent} we obtain:
\begin{lemma} \label{lem:gf:1n:basic}
	Consider the notation in \cref{subsec:1n:notation}. Then
	\begin{enumerate} [label = (\roman*)]
		\item \label{lem:gf:1n:basic:item1}
		it holds that $\Theta \in C^1 ( [0, \infty ) , \R^3)$, 
		\item \label{lem:gf:1n:basic:item2}
		it holds for all $t \in [0 , \infty )$ that $\Theta ( t ) \in \cM $,
		and
		\item \label{lem:gf:1n:basic:item3}
		it holds that $[0 , \infty ) \ni t \mapsto \cL ( \Theta ( t ) ) \in \R$ is non-increasing.
	\end{enumerate}
\end{lemma}

\begin{remark}
	It is not hard to see that for all $\theta \in \cM$ the following properties hold:
	\begin{itemize}
		\item If $\theta_1 > 0$ then $I^\theta = (q^\theta , \infty ) \cap [0 , 1 ]$,
		and if $\theta_1 < 0$ then $I^\theta = ( - \infty , q^\theta ) \cap [0 , 1 ]$.
		\item It holds that $\mu ( I^\theta ) \in (0 , 1 )$ if and only if $q^\theta \in (0 , 1 )$.
		\item It holds that $\mu ( I^\theta ) = 0$ if and only if $I^\theta = \emptyset$.
		\item It holds that $\mu ( I^\theta ) = 1$ if and only if $I^\theta \supseteq (0 , 1 )$.
	\end{itemize}
	This and the fact that $\Theta \in C ( [ 0 , \infty ) , \R^3 )$ easily imply that $[0 , \infty ) \ni t \mapsto \mu ( I ^{ \Theta ( t ) } ) \in \R$ is continuous.
\end{remark}

\subsection{Proof of the boundedness in simple cases}
We
first show the following:

\begin{lemma} \label{lem:bounded:case:simple1}
	Consider the notation in \cref{subsec:1n:notation}.
	Then for every $\varepsilon > 0$ it holds that
	\begin{equation}
	\sup\nolimits_{t \in [0 , \infty )} \br*{ \norm{\Theta(t ) }  \indicator{ [\varepsilon ,  1 ) } ( \mu ( I ^{ \Theta ( t ) } ) )  } < \infty .
	\end{equation}
\end{lemma}

Notice that, due to \cref{lem:gf:1n:basic}, it suffices to bound $\Theta_3 ( t ) $.
For this we use the following elementary lemma,
which is proved, e.g., in \cite[Corollary 5.2]{JentzenRiekertFlow}.

\begin{lemma} \label{lem:affine:integral}
	Let $\alpha, \beta \in \R$
	and let $I \subseteq \R$ be a bounded interval.
	Then $\int_I ( \alpha x + \beta ) ^2 \, \d x \ge \tfrac{\alpha^2 }{12} ( \mu ( I ) ) ^3$.
\end{lemma}

\begin{cproof}{lem:bounded:case:simple1}
	Throughout this proof let $\varepsilon > 0$,
	let $\fC = \norm{ f - \overline{f} } _{L^2 ( [0,1] ) }$,
	and let $\cT \subseteq[0 , \infty )$
	satisfy
	$\cT = \cu{t \in [0 , \infty ) \colon \mu ( I ^{ \Theta ( t ) } ) \in [ \varepsilon , 1 ) } $.
	\Nobs that for all $t \in \cT$ it holds that
	$q^{ \Theta ( t ) } \in (0 , 1 )$
	and
	$\abs{\Theta_2 ( t ) } \le \abs{\Theta_1 ( t ) }$. Hence, we obtain for all $t \in \cT$ that $\abs{\Theta_1 ( t ) } \ge 2^{- \nicefrac{1}{2} }$.
	Furthermore, the triangle inequality \proves for all $\theta \in \R^3 $
	that
	\begin{equation}
	\begin{split}
	\sqrt{\cL ( \theta ) } 
	&\ge \rbr*{ \int_0^1  \rbr[\big]{ \theta_3 ( \max \cu{\theta_1 s + \theta_2 , 0 } - m ( \theta ) ) } ^2 \, \d s } ^{ \nicefrac{1}{2} } - \fC \\
	&\ge \rbr*{ \abs{\theta_3 } ^2 \int_{I^\theta } ( \theta_1 s + \theta_2 - m ( \theta ) ) ^2 \, \d s } ^{ \nicefrac{1}{2 } } - \fC .
	\end{split}
	\end{equation}
	Combining this with \cref{lem:affine:integral},
	\cref{lem:gf:1n:basic:item3} in \cref{lem:gf:1n:basic},
	and the fact that $\forall \, t \in \cT \colon \abs{\Theta_1 ( t ) } \ge 2^{- \nicefrac{1}{2} }$
	\proves for all $t \in \cT$
	that
	\begin{equation}
	\begin{split}
	\sqrt{\cL (\Theta ( 0 )  ) }
	& \ge \sqrt{ \cL ( \Theta ( t ) ) }
	\ge \abs{\Theta_3 ( t ) } \rbr*{ \int_{I^{\Theta ( t ) } } \rbr[\big]{ \Theta_1 ( t )  s + \Theta_2 ( t )  - m ( \Theta ( t ) ) } ^2 \, \d s } ^{ \nicefrac{1}{2 } } - \fC \\
	& \ge \abs{\Theta_3 ( t ) } \abs{\Theta_1 ( t ) } ( \tfrac{1}{12} \mu ( I^{\Theta ( t ) } ) ^3 ) ^{ \nicefrac{1}{2}} - \fC \ge 24^{ - \nicefrac{1}{2} } \abs{\Theta_3 ( t ) } \varepsilon ^{ \nicefrac{3}{2}} - \fC .
	\end{split}
	\end{equation}
	This \proves that $\sup_{t \in \cT } \abs{\Theta_3 ( t ) } < \infty $.
\end{cproof}

From \cref{lem:bounded:case:simple1} we obtain the boundedness if $I ^{ \Theta ( t ) }$ is not the entire interval $(0,1)$, but has a positive measure bounded away from zero.

\begin{lemma}  \label{lem:bounded:case:simple2}
	Consider the notation in \cref{subsec:1n:notation}.
	Then for every $\varepsilon > 0$ it holds that
	\begin{equation}
	\sup\nolimits_{t \in [0 , \infty ) } \br[\big]{ \norm{\Theta ( t ) } \indicator{[\varepsilon , \infty ) } ( \abs{\Theta_1 ( t ) } ) \indicator{\cu{1} } ( \mu ( I^{ \Theta ( t ) } ) ) } < \infty 
	\end{equation}
\end{lemma}

\begin{cproof}{lem:bounded:case:simple2}
	Throughout this proof let $\varepsilon > 0$,
	let $\fC = \norm{ f - \overline{f} } _{L^2 ( [0,1] ) }$,
	and let $\cT \subseteq[0 , \infty )$
	satisfy
	$\cT = \cu{t \in [0 , \infty ) \colon \abs{\Theta_1 ( t ) } \ge \varepsilon , \, \mu ( I^{ \Theta ( t ) } ) = 1  } $.
	Using the same arguments as in the proof of \cref{lem:bounded:case:simple2} yields for all $t \in \cT$ that
	\begin{equation}
	\begin{split}
	\sqrt{\cL (\Theta ( 0 ) ) }
	& \ge \sqrt{ \cL ( \Theta ( t ) ) }
	\ge \abs{\Theta_3 ( t ) } \rbr*{ \int_0^1  \rbr[\big]{ \Theta_1 ( t )  s + \Theta_2 ( t )  - m ( \Theta ( t ) ) } ^2 \, \d s  } ^{ \nicefrac{1}{2 } } - \fC \\
	& \ge ( \tfrac{1}{12}  ) ^{ \nicefrac{1}{2}} \abs{\Theta_3 ( t ) } \abs{\Theta_1 ( t ) }  - \fC \ge \tfrac{1}{4} \abs{\Theta_3 ( t ) } \varepsilon - \fC .
	\end{split}
	\end{equation}
	This \proves that $\sup_{t \in \cT } \abs{\Theta_3 ( t ) } < \infty $.
\end{cproof}

\begin{prop}
	\label{prop:neuron:active:bounded}
	Consider the notation in \cref{subsec:1n:notation}.
	Then it holds for all $t \in [0 , \infty )$ with $q^{ \Theta ( t ) } \notin ( 0 , 1 )$
	and $\abs{\Theta_1 ( t ) } < 1 $
	that
	\begin{equation}
	\label{prop:neuron:active:bounded:eqclaim}
	\frac{\d}{\d t } \rbr*{ \abs{\Theta_3 ( t ) } ^2 + \ln ( 1 - \abs { \Theta_1 ( t ) } ^2 ) } = 0.
	\end{equation}
\end{prop}

\begin{cproof}{prop:neuron:active:bounded}
	\Nobs that for all $t \in [0 , \infty )$ with $\abs{\Theta_1 ( t ) } < 1 $ we have that
	\begin{equation} 
	\label{prop:neuron:active:bounded:eq1}
	\begin{split}
	\frac{\d}{\d t } \rbr*{ \abs{\Theta_3 ( t ) } ^2 + \ln ( 1 - \abs { \Theta_1 ( t ) } ^2 ) } 
	&= 2  \Theta_3 ( t ) \fG_3 ( \Theta ( t ) ) - 2 \frac{\Theta_1 ( t ) }{1 - \abs{\Theta_1 ( t ) } ^2 } \fG_1 ( \Theta ( t ) ) \\
	& = 2  \Theta_3 ( t ) \fG_3 ( \Theta ( t ) ) - 2 \frac{\Theta_1 ( t ) }{ \abs{\Theta_ 2 ( t ) } ^2 } \fG_1 ( \Theta ( t ) ).
	\end{split}
	\end{equation}
	Furthermore, if $q^{ \Theta ( t ) } \notin ( 0 , 1 )$ we either have $I^{ \Theta ( t ) } = \emptyset$
	or $I^{ \Theta ( t ) }  \supseteq (0 , 1 ) $.
	In the first case \cref{lem:gradient:basic} \proves that $\fG_1 ( \Theta ( t ) ) = \fG_3 ( \Theta ( t ) ) = 0$.
	 \cref{prop:neuron:active:bounded:eq1} therefore \proves \cref{prop:neuron:active:bounded:eqclaim}.
	
	In the second case we obtain from \cref{lem:gradient:basic} that
	\begin{equation}
	\fG_1 ( \Theta ( t ) ) = 2 ( \Theta_2 ( t ) ) ^2 \Theta_3 ( t ) \int_0^1 \rbr*{ \Theta_3 ( t ) ( \max \cu{\Theta_1 ( t ) s + \Theta_2 ( t ) , 0 } - m ( \Theta ( t ) ) ) + \overline{f} - f ( s )   } s \, \d s
	\end{equation}
	and
	\begin{equation}
	\fG_3 ( \Theta ( t ) ) =  2 \Theta_1 ( t ) \Theta_3 ( t ) \int_0^1 \rbr*{ \Theta_3 ( t ) ( \max \cu{\Theta_1 ( t ) s + \Theta_2 ( t ) , 0 } - m ( \Theta ( t ) ) ) + \overline{f} - f ( s )   } s \, \d s
	\end{equation}
	Combining this with \cref{prop:neuron:active:bounded:eq1} establishes \cref{prop:neuron:active:bounded:eqclaim}.
\end{cproof}

\begin{remark}
	An analogous statement to \cref{prop:neuron:active:bounded} can be proved for any number $H \in \N$ of neurons on the hidden layer, using similar identities for the gradient components.
\end{remark} 

Using the last two results, we get boundedness in the case $q^{ \Theta ( t ) } \notin (0 , 1 )$. Indeed, if $I^{ \Theta ( t ) } = \emptyset$ then $\fG ( \Theta ( t ) ) = 0$, so it cannot diverge.
If $I ^{ \Theta ( t ) } \supseteq (0 , 1 )$ and $\abs{\Theta_1 ( t ) } $ is bounded away from zero the boundedness follows from \cref{lem:bounded:case:simple2}.
If 
$I ^{ \Theta ( t ) } \supseteq (0 , 1 )$ and $\abs{\Theta_1 ( t ) } $ is bounded away from $1$ the boundedness follows from
\cref{prop:neuron:active:bounded}.

The remaining and more difficult cases occur when $I^{ \Theta ( t ) }$ has small positive measure. This is the content of the next two subsubsections.

\subsection{The case that the breakpoint is close to 1}

In this subsection we will deal with the case that the activity interval
 $I^{ \Theta ( t ) }$ is non-empty and contained in some interval $[1 - \varepsilon , 1 ]$ for a small $\varepsilon > 0$, which is not covered by the previous results.
\Nobs that $I^{ \Theta ( t ) } $ can only be of the considered form if $q ^{ \Theta ( t ) } \in (0 , 1 ) $ and $ \Theta_2 ( t ) < 0 < \Theta_1 ( t )  $.
Furthermore, we have $\abs{\Theta_2 ( t ) } ^2 < \frac{1}{2} < \abs{\Theta_1 ( t ) } ^2 $. This will be used throughout this section.

\begin{lemma} \label{lem:integral:formulae}
	Consider the notation in \cref{subsec:1n:notation} and let $\theta \in \cM $ satisfy $q^\theta \in (0 , 1 )$ and $\theta_1 > 0$.
	Then
	\begin{enumerate} [ label = (\roman*) ]
		\item
		\label{lem:integral:formulae:item1}
		it holds that $m ( \theta ) = \frac{\theta_1}{2} ( 1 - q^\theta ) ^2 $,
		\item
		\label{lem:integral:formulae:item2}
		it holds that 
		\begin{equation}
		\int_{I^\theta } \rbr*{ \max \cu{\theta_1 s + \theta_2 , 0 } - m ( \theta ) } \, \d s = \tfrac{\theta_1}{2} ( 1 - q^\theta ) ^2 q ^\theta ,
		\end{equation}
		\item
		\label{lem:integral:formulae:item3}
		it holds that
		\begin{equation}
		\int_0^1 \rbr*{ \max \cu{\theta_1 s + \theta_2 , 0 } - m ( \theta ) } ^2 \, \d s = \theta_1 ^2 ( 1 - q^\theta ) ^3 ( \tfrac{1}{12} + \tfrac{q^\theta }{4} ) ,
		\end{equation}
		and
		\item 
		\label{lem:integral:formulae:item4}
		it holds that
		\begin{equation}
		\begin{split}
		\fG_1 ( \theta ) &=
		2 \theta_3 \rbr*{ \frac{\theta_1 \theta_2^2 \theta_3}{12} ( 1 - q^\theta ) ^2 ( 7 + 2 q^\theta + 3 ( q^\theta ) ^2 )
			+ \int_{q^\theta} ^1 ( \overline{f} - f ( s ) ) ( \theta_2^2 s - \theta_1 \theta_2 ) \, \d s } , \\
		\fG_3 ( \theta )
		&=
		2 \theta_1^2 \theta_3 ( 1 - q^\theta ) ^3 \rbr*{ \frac{1}{12} + \frac{q^\theta }{4} } + 2 \int_{q^\theta } ^1 ( \overline{f} - f ( s ) ) \max \cu{\theta_1 s + \theta_2 , 0 } \, \d s  .
		\end{split}
		\end{equation}
	\end{enumerate}
\end{lemma}

\begin{cproof}{lem:integral:formulae}
	First, we have
	\begin{equation}
	\begin{split}
	m ( \theta ) 
	&= \int_0^1 \max \cu{\theta_1 s + \theta_2 , 0 } \, \d s
	= \theta_1 \int_{q^\theta} ^1 ( s - q^\theta ) \, \d s 
	= \tfrac{\theta_1}{2} ( 1 - q^\theta ) ^2.
	\end{split}
	\end{equation}
	This establishes \cref{lem:integral:formulae:item1}.
	Next, \cref{lem:integral:formulae:item1} implies that
	\begin{equation}
	\begin{split}
	\int_{I^\theta } (\max \cu{\theta_1 s + \theta_2 , 0 } - m ( \theta ) ) \, \d s 
	&= m ( \theta ) - \int_{q^\theta} ^1 m ( \theta) \, \d s
	= m ( \theta ) - (1 - q^\theta ) m ( \theta ) \\
	&= q^\theta m ( \theta ) = \tfrac{\theta_1}{2} ( 1 - q^\theta ) ^2 q ^\theta.
	\end{split}
	\end{equation}
	This establishes \cref{lem:integral:formulae:item2}.
	Moreover,
	\nobs that
	\begin{equation}
	\begin{split}
	& \int_0^1 ( \max \cu{\theta_1 s + \theta_2 , 0 } - m ( \theta ) ) ^2 \, \d s \\
	& = \int_0^1 ( \max \cu{\theta_1 s + \theta_2 , 0 } ) ^2 \, \d s
	- 2 m ( \theta ) \int_0^1 \max \cu{\theta_1 s + \theta_2 , 0 } \, \d s 
	+ \int_0^1 m ( \theta ) ^2 \, \d s \\
	&= \theta_1^2 \int_{q^\theta }^1 ( s - q^\theta ) ^2 \, \d s  - 2 m ( \theta ) ^2 + m ( \theta ) ^2 
	= \theta_1^2 \int_0^{ 1 - q^\theta } s^2 \, \d s - m ( \theta ) ^2 
	\\
	&= \tfrac{\theta_1^2}{3} ( 1 - q^\theta ) ^3 - \tfrac{\theta_1^2}{4} ( 1 - q^\theta ) ^4 
	= \theta_1 ^2 ( 1 - q^\theta ) ^3 ( \tfrac{1}{3} - \tfrac{ 1 - q^\theta }{4} )
	= \theta_1 ^2 ( 1 - q^\theta ) ^3 ( \tfrac{1}{12} + \tfrac{q^\theta }{4} ) .
	\end{split}
	\end{equation}
	This establishes \cref{lem:integral:formulae:item3}.
	In addition, \cref{eq:gen:gradient:explicit} assures that
	\begin{equation}
	\begin{split}
	\fG_3 ( \theta ) &= 2 \theta_3 \int_0^1 ( \max \cu{\theta_1 s + \theta_2 , 0 } - m ( \theta ) ) ^2 \, \d s \\
	& \quad + 2 \int_{ 0 } ^1 ( \overline{f} - f ( s ) ) ( \max \cu{\theta_1 s + \theta_2 , 0 } - m ( \theta ) ) \, \d s \\
	&= 2 \theta_1^2 \theta_3 ( 1 - q^\theta ) ^3 \rbr*{ \frac{1}{12} + \frac{q^\theta }{4} } + 2 \int_{q^\theta } ^1 ( \overline{f} - f ( s ) ) \max \cu{\theta_1 s + \theta_2 , 0 } \, \d s .
	\end{split}
	\end{equation}
	Furthermore, \nobs that
	\begin{equation}
	\begin{split}
	& \int_{I^\theta } ( \max \cu{\theta_1 s + \theta_2 , 0 } - m ( \theta ) ) ( \theta_2^2 s - \theta_1 \theta_2 ) \, \d s \\
	&= - \theta_1 \theta_2 \int_{I^\theta } ( \max \cu{\theta_1 s + \theta_2 , 0 } - m ( \theta ) ) \, \d s + \theta_2^2 \int_{I^\theta} s ( \theta_1 s + \theta_2 ) \, \d s - \theta_2^2 m ( \theta ) \int_{I^\theta } s \, \d s \\
	&= - \tfrac{\theta_1^2 \theta_2}{2} ( 1 - q^\theta ) ^2 q^\theta + \theta_2^2 \theta_1 \rbr*{ \tfrac{1 - (q^\theta ) ^3}{3} - \tfrac{q^\theta ( 1 - (q^\theta ) ^2 ) }{2} } - \tfrac{\theta_2^2 \theta_1}{2} ( 1 - q^\theta ) ^2 \rbr*{ \tfrac{1 - (q^\theta ) ^2 }{2} } \\
	&= ( 1 - q^\theta ) \theta_1 \theta_2^2
	\rbr*{\tfrac{1 - q^\theta }{2} + \tfrac{1 + q^\theta + ( q^\theta ) ^2 }{3} - \tfrac{q^\theta ( 1 + q^\theta ) }{2} - \tfrac{ ( 1 - q^\theta ) ( 1-(q^\theta ) ^2 )}{4} }
	\\
	&= (1 - q^\theta ) ^2 \tfrac{\theta_1 \theta_2^2}{12} ( 7 + 2 q^\theta + 3 ( q^\theta ) ^2 ).
	\end{split}
	\end{equation}
	Combining this with \cref{eq:gen:gradient:explicit} establishes \cref{lem:integral:formulae:item4}.
\end{cproof}

Next, by symmetry we may assume wlog that $\overline{f} \le f ( 1 ) $. (Otherwise replace $f \with -f$, $\theta_3 \with - \theta_3 $.)

\begin{lemma} \label{lem:endpoint:mean:value}
	Consider the notation in \cref{subsec:1n:notation},
	assume $\overline{f} = f ( 1 )$,
	and assume that $f$ is Lipschitz continuous.
	Then there exists $c \in \R$
	which satisfy for all $t \in [0, \infty )$ with $  I^{\Theta ( t ) } \subseteq [ \frac{1}{2} , 1 ]$
	and $ \abs{ \Theta_3 ( t ) } \ge c $
	that
	\begin{equation}
	\frac{\d }{\d t } \abs{\Theta_3 ( t ) } ^2  \le 0.
	\end{equation}
\end{lemma}

\begin{cproof}{lem:endpoint:mean:value}
	First, the assumption that $f$ is Lipschitz continuous ensures that there exists $L \in (0 , \infty )$
	which satisfies for all $s \in [0 , 1 ]$ that $\abs{ f ( s ) - f ( 1 ) } = \abs{f( s ) - \overline{f} } \le L ( 1 - s )$.
	Combining this with \cref{lem:integral:formulae}
	demonstrates for all $\theta \in \cM $
	with $ \emptyset \not= I^{ \theta } \subseteq [ \frac{1}{2} , 1 ]$
	that
	\begin{equation}
	\begin{split}
	\abs*{ \int_{q^\theta } ^1 ( \overline{f} - f ( s ) ) \max \cu{\theta_1 s + \theta_2 , 0 } \, \d s }
	\le L \theta_1 \int_{q^\theta}^1 ( 1 - s ) ( s - q^\theta ) \, \d s
	= \tfrac{L \theta_1}{6} ( 1 - q^\theta ) ^3 .
	\end{split}
	\end{equation}
	This,
	the chain rule,
	and the fact that for all $t \in [0 , \infty )$ with $\emptyset \not= I ^{ \Theta ( t ) } \subseteq [ \frac{1}{2} , 1 ]$ it holds that $2^{ - \nicefrac{1}{2}} \le \abs{\Theta_1 ( t ) } \le 1$ show that for all $t \in [0 , \infty )$ with $\emptyset \not= I ^{ \Theta ( t ) } \subseteq [ \frac{1}{2} , 1 ]$
	we have that
	\begin{equation}
	\begin{split}
	\frac{\d }{\d t } \abs{\Theta_3 ( t ) } ^2 
	&= - 2 \Theta_3 ( t ) \fG _3 ( \Theta ( t ) ) \\
	&\le - 2 \abs{\Theta_1 ( t ) \Theta_3 ( t ) } ^2 ( 1 - q^{ \Theta ( t ) } ) ^3  \rbr*{ \tfrac{1}{12} + \tfrac{q^{\Theta ( t ) } }{4} }
	+ L \tfrac{ \abs{\Theta_1 ( t ) \Theta_3 ( t) } }{3}  ( 1 - q^{\Theta ( t ) } ) ^3   \\
	& \le ( 1 - q^{\Theta ( t ) } ) ^3 \rbr*{ - \tfrac{\abs{\Theta_3 ( t ) } ^2 }{ 12 } + \tfrac{L \abs{\Theta_3 ( t ) } }{3 } }
	= ( 1 - q^{\Theta ( t ) } ) ^3 \tfrac{\abs{\Theta_3 ( t ) } } { 12 } ( 4 L - \abs{\Theta_3 ( t ) } ) .
	\end{split}
	\end{equation}
	Hence, we obtain for all $t \in [0 , \infty )$ with $I ^{ \Theta ( t ) } \subseteq [ \frac{1}{2} , 1 ]$ and $\abs{\Theta_3 ( t ) } \ge 4 L$ that
	$\frac{\d }{\d t } \abs{\Theta_3 ( t ) } ^2  \le 0$.
\end{cproof}

\begin{lemma}
	\label{lem:difficult:case1}
	Consider the notation in \cref{subsec:1n:notation} and assume $\overline{f} < f ( 1 )$.
	Then there exists $\varepsilon  \in (0 , \nicefrac{1}{2} )$
	which satisfies for all $t \in [0, \infty )$ with $ \emptyset \not= I^{\Theta ( t ) } \subseteq [ 1 - \varepsilon , 1 ]$ and $\Theta_3 (t ) \le 0$
	that
	\begin{equation}
	\frac{\d }{\d t } \abs{\Theta_3 ( t ) } ^2  \le 0.
	\end{equation}
\end{lemma}

\begin{cproof}{lem:difficult:case1}
	\Nobs that the fact that $f$ is continuous assures that there exists $\varepsilon  \in (0 , \nicefrac{1}{2} )$
	which satisfies for all $s \in [ 1 - \varepsilon , 1 ]$ that $f(s ) > \overline{f}$.
	This implies for all $t \in [0, \infty )$ with $ \emptyset \not= I^{\Theta ( t ) } \subseteq [ 1 - \varepsilon , 1 ]$
	that $\int_{q^{\Theta ( t ) } }^1 ( \overline{f} - f ( s ) ) \max \cu{\Theta_1 ( t ) s + \Theta_2 ( t ) , 0 } \, \d s \le 0$.
	Combining this with \cref{lem:integral:formulae} demonstrates for all $t \in [0, \infty )$ with $ \emptyset \not= I^{\Theta ( t ) } \subseteq [ 1 - \varepsilon , 1 ]$ and $\Theta_3 (t ) \le 0$
	that
	$\frac{\d }{\d t } \abs{\Theta_3 ( t ) } ^2 
	= - 2 \Theta_3 ( t ) \fG _3 ( \Theta ( t ) ) \le 0$.
\end{cproof}

\begin{lemma}
	\label{lem:difficult:case2}
	Consider the notation in \cref{subsec:1n:notation} and assume $\overline{f} < f ( 1 )$.
	Then there exists $\varepsilon  \in (0 , \nicefrac{1}{2} )$
	which satisfies for all $t \in [0, \infty )$ with $ \emptyset \not= I^{\Theta ( t ) } \subseteq [ 1 - \varepsilon , 1 ]$ and $\Theta_3 ( t ) > 0$
	that
	\begin{equation}
	\frac{\d }{\d t } \rbr*{ \abs{\Theta_3 ( t ) } ^2 - \tfrac{5}{8 } \abs{ \Theta_1 ( t ) - 2^{- \nicefrac{1}{2} } } ^2 } \le 0.
	\end{equation}
\end{lemma}

\begin{cproof}{lem:difficult:case2}
	First, the fact that $f$ is continuous ensures that there exist $\beta \in (0 , \infty )$,
	$\varepsilon  \in (0 , \nicefrac{1}{2} )$
	which satisfy for all $s \in [ 1- \varepsilon , 1 ]$
	that
	$\overline{f} + \beta < f ( s ) < \overline{f} + \frac{10 \beta }{9 }$.
	This implies for all $\theta \in \cM $ with $\emptyset \not= I^\theta \subseteq [ 1 - \varepsilon , 1 ]$
	that
	\begin{equation}
	\int_{I^\theta } ( f ( s ) - \overline{f} ) \max \cu{\theta_1 s + \theta_2 , 0 } \, \d s \le \frac{10 \beta}{9} \int_{q^\theta} ^1  \max \cu{\theta_1 s + \theta_2 , 0 } \, \d s = \frac{5 \beta}{9} \theta_1 ( 1 - q^\theta ) ^2 
	\end{equation}
	and
	\begin{equation}
	\begin{split}
	\int_{I^\theta} ( f ( s ) - \overline{f} ) ( \theta_2^2 s - \theta_1 \theta_2 ) \, \d s 
	\ge \frac{\beta}{2} \theta_1^2 q^\theta ( 1 - q^\theta ) ( 2 + q^\theta + (q^{\theta} ) ^2 ) .
	\end{split}
	\end{equation}
	Combining this with \cref{lem:integral:formulae}
	demonstrates for all
	$\theta \in \cM $ with $ \emptyset \not= I^\theta \subseteq [ 1 - \varepsilon , 1 ]$
	and
	$\theta_3 \ge 0$
	that
	\begin{equation}
	\fG_3 ( \theta ) \ge \theta_1 ( 1 - q^\theta ) ^2 \rbr*{ 2 \theta_1 \theta_3 ( 1 - q^\theta ) \rbr*{ \tfrac{1}{12} + \tfrac{q^\theta}{4} } - \tfrac{10 \beta}{9} }
	\end{equation}
	and
	\begin{equation}
	\begin{split}
	\fG_1 ( \theta ) 
	&\le 2 \theta_3 \rbr*{ \theta_1 \theta_2^2 \theta_3 ( 1 - q^\theta ) ^2 - \tfrac{\beta}{2} \theta_1^2 q^\theta ( 1 - q^\theta ) ( 2 + q^\theta + (q^{\theta} ) ^2 ) } \\
	& \le \theta_1 \theta_3  ( 1 - q^\theta ) \rbr*{\theta_3 ( 1 - q^\theta ) - \beta \theta_1 q^\theta  ( 2 + q^\theta + ( q^\theta ) ^2 ) }.
	\end{split}
	\end{equation}
	In addition,
	\nobs that for all $\theta \in \cM $ with $ \emptyset \not= I^\theta \subseteq [ 1 - \varepsilon , 1 ]$ it holds that
	\begin{equation}
	\begin{split}
	\theta_1 - 2^{- \nicefrac{1}{2} } 
	&= \frac{\theta_1}{\sqrt{2}} \rbr*{ \sqrt{2} - \theta_1^{-1} \sqrt{ \theta_1^2 + \theta_2^2 } }
	= \frac{\theta_1}{\sqrt{2} } \rbr*{ \sqrt{2} - \sqrt{1 + (q^\theta ) ^2 } } \\
	&= \frac{\theta_1 ( 1 - ( q^\theta ) ^2  ) }{2 + \sqrt{2 ( 1 +(q^\theta ) ^2) } } .
	\end{split}
	\end{equation}
	This and the chain rule show for all 
	$t \in [0 , \infty )$
	with $ \emptyset \not= I^{\Theta ( t ) } \subseteq [ 1 - \varepsilon , 1 ]$ and $\Theta_3 ( t ) > 0$
	that
	\begin{equation}
	\begin{split}
	&  \frac{\d }{\d t } \rbr*{ \abs{\Theta_3 ( t ) } ^2 - \tfrac{5}{8} \abs{ \Theta_1 ( t ) - 2^{- \nicefrac{1}{2} } } ^2 } \\
	&= - 2 \Theta_3 ( t ) \fG_3 ( \Theta ( t ) )
	+ \tfrac{5}{4} \fG_1 ( \Theta ( t ) ) ( \Theta_1 ( t ) - 2^{- \nicefrac{1}{2} } )  \\
	& \le 
	2 \Theta_1 ( t ) \Theta_3 ( t ) ( 1 - q^{ \Theta ( t ) } ) ^2 
	\Big( \tfrac{10 \beta }{9} - \Theta_1 ( t ) \Theta_3 ( t ) ( 1 - q^{ \Theta ( t ) } ) \rbr*{ \tfrac{1}{6} + \tfrac{q^{ \Theta ( t ) } }{2} } \\
	& \qquad - \tfrac{ 5 \Theta_1 ( t ) ( 1 + q^{ \Theta ( t ) } ) } { 8 + 4 \sqrt{2 ( 1 +(q^{ \Theta ( t ) } ) ^2) } } \rbr*{ \beta \Theta_1 ( t ) q^{ \Theta ( t ) } ( 2 + q^{ \Theta ( t ) } + ( q ^{ \Theta ( t ) } ) ^2 ) - \Theta_3 ( t ) ( 1 - q^{ \Theta ( t ) } )  } \Big) .
	\end{split}
	\end{equation}
	Next \nobs that the fact that $\forall \, q \in [0 , 1 ] \colon \frac{1+q}{2 + \sqrt{2 ( 1+q^2 ) } } \le \frac{1}{2}$
	ensures that there exists $\eta \in (0 , \varepsilon )$
	which satisfies for all
	$t \in [0 , \infty )$
	with $ \emptyset \not= I^{\Theta ( t ) } \subseteq [ 1 - \eta  , 1 ]$
	that
	\begin{multline}
	\tfrac{ 5 \abs{ \Theta_1 ( t ) }^2 ( 1 + q^{ \Theta ( t ) } ) } { 8 + 4 \sqrt{2 ( 1 +(q^{ \Theta ( t ) } ) ^2) } }
	q^{ \Theta ( t ) } ( 2 + q^{ \Theta ( t ) } + ( q ^{ \Theta ( t ) } ) ^2 ) \ge \tfrac{10}{9} , \\ 
	\tfrac{ 5 ( 1 + q^{ \Theta ( t ) } ) } { 8 + 4 \sqrt{2 ( 1 +(q^{ \Theta ( t ) } ) ^2) } } \le \tfrac{5}{8},
	\qqandqq \tfrac{1}{6} + \tfrac{q^{ \Theta ( t ) } }{2} > \tfrac{5}{8} .
	\end{multline}
	Therefore,
	we obtain for all $t \in [0 , \infty )$
	with $ \emptyset \not= I^{\Theta ( t ) } \subseteq [ 1 - \eta  , 1 ]$ and $\Theta_3 ( t ) > 0$
	that
	\begin{equation}
	\begin{split}
	&  \frac{\d }{\d t } \rbr*{ \abs{\Theta_3 ( t ) } ^2 - \tfrac{5}{8} \abs{ \Theta_1 ( t ) - 2^{- \nicefrac{1}{2} } } ^2 } \\
	& \le  2 \Theta_1 ( t ) \Theta_3 ( t ) ( 1 - q^{ \Theta ( t ) } ) ^2 \Big(  \tfrac{10 \beta }{9} -  \Theta_1 ( t ) \Theta_3 ( t ) ( 1 - q^{ \Theta ( t ) } ) \rbr*{ \tfrac{1}{6} + \tfrac{q^{ \Theta ( t ) } }{2} } \\
	& \qquad - \tfrac{10 \beta }{9} + \tfrac{5}{8} \Theta_1 ( t ) \Theta_3 ( t ) ( 1 - q^{ \Theta ( t ) } ) \Big) \le  0.
	\end{split}
	\end{equation}
\end{cproof}

\subsection{The case that the breakpoint is close to 0}
Finally, we consider the case where the activity interval $I^{ \Theta ( t ) }$ is non-empty and contained in some interval $[0 , \varepsilon ]$ with $\varepsilon > 0$ small.
The arguments are essentially analogous to the previous case. Note that this time we must have $q ^{ \Theta ( t ) } \in (0 , 1 )$ and $\Theta_2  ( t ) > 0 > \Theta_1 ( t ) $. Furthermore, for small $\varepsilon >0$ we have that $\Theta_1 ( t ) $ is close to $-1$ and $\Theta_2 ( t )$ is close to $0$.

\begin{lemma} \label{lem:integral:formulae2}
	Consider the notation in \cref{subsec:1n:notation} 
	and let $\theta \in \cM $ satisfy $q^\theta \in (0 , 1 )$ and $\theta_1 < 0$.
	Then
	\begin{enumerate} [ label = (\roman*) ]
		\item
		\label{lem:integral:formulae2:item1}
		it holds that $m ( \theta ) = \frac{\theta_1}{2} ( q^\theta ) ^2 $,
		\item
		\label{lem:integral:formulae2:item2}
		it holds that 
		\begin{equation}
		\int_{I^\theta } (\max \cu{\theta_1 s + \theta_2 , 0 } - m ( \theta ) ) \, \d s = \tfrac{\theta_1}{2} ( 1 - q^\theta ) ( q ^\theta ) ^2 ,
		\end{equation}
		\item
		\label{lem:integral:formulae2:item3}
		it holds that
		\begin{equation}
		\int_0^1 ( \max \cu{\theta_1 s + \theta_2 , 0 } - m ( \theta ) ) ^2 \, \d s = \theta_1 ^2 ( q^\theta ) ^3 ( \tfrac{1}{3} - \tfrac{q^\theta }{4} ) ,
		\end{equation}
		and
		\item 
		\label{lem:integral:formulae2:item4}
		it holds that
		\begin{equation}
		\begin{split}
		\fG_1 ( \theta ) &=
		2 \theta_3 \rbr*{ -  \frac{\theta_1^3 \theta_3}{12} ( q^\theta ) ^3 ( 6 + 6 q^\theta + 2 ( q^\theta ) ^2 + 3 ( q^\theta ) ^ 3 )
			+ \int_0^{q^\theta}  ( \overline{f} - f ( s ) ) ( \theta_2^2 s - \theta_1 \theta_2 ) \, \d s } , \\
		\fG_3 ( \theta )
		&=
		2 \theta_1^2 \theta_3 (   q^\theta ) ^3 \rbr*{ \frac{1}{3} - \frac{q^\theta }{4} } + 2 \int_0^{q^\theta } ( \overline{f} - f ( s ) ) \max \cu{\theta_1 s + \theta_2 , 0 } \, \d s .
		\end{split}
		\end{equation}
	\end{enumerate}
\end{lemma}

In the following consider the case $\overline{f} \le f( 0 )$, the case $\overline{f} > f( 0 )$ being analogous.

\begin{lemma} \label{lem:endpoint:mean:value2}
	Consider the notation in \cref{subsec:1n:notation}, assume $\overline{f} = f ( 0 )$,
	and assume that $f$ is Lipschitz continuous.
	Then there exists $c \in \R$
	which satisfy for all $t \in [0, \infty )$ with $  I^{\Theta ( t ) } \subseteq [ 0 , \frac{1}{2} ]$
	and $ \abs{ \Theta_3 ( t ) } \ge c $
	that
	\begin{equation}
	\frac{\d }{\d t } \abs{\Theta_3 ( t ) } ^2  \le 0.
	\end{equation}
\end{lemma}

\begin{cproof}{lem:endpoint:mean:value2}
	First, the assumption that $f$ is Lipschitz continuous ensures that there exists $L \in (0 , \infty )$
	which satisfies for all $s \in [0 , 1 ]$ that $\abs{ f ( s ) - f ( 0 ) } = \abs{f( s ) - \overline{f} } \le L  s$.
	Combining this with \cref{lem:integral:formulae2}
	demonstrates for all $\theta \in \cM $
	with $ \emptyset \not= I^{ \theta } \subseteq [ 0 , \frac{1}{2} ]$
	that
	\begin{equation}
	\begin{split}
	\abs*{ \int_0^{q^\theta }  ( \overline{f} - f ( s ) ) \max \cu{\theta_1 s + \theta_2 , 0 } \, \d s }
	\le L \abs{ \theta_1 } \int_0^{q^\theta}  s  (   q^\theta - s  ) \, \d s
	= \tfrac{L \theta_1}{6} (  q^\theta ) ^3 .
	\end{split}
	\end{equation}
	This and the chain rule show for all $t \in [0 , \infty )$ with $\emptyset \not= I ^{ \Theta ( t ) } \subseteq [ 0 , \frac{1}{2}  ]$
	that
	\begin{equation}
	\begin{split}
	\frac{\d }{\d t } \abs{\Theta_3 ( t ) } ^2 
	&= - 2 \Theta_3 ( t ) \fG _3 ( \Theta ( t ) ) \\
	&\le - 2 \abs{\Theta_1 ( t ) \Theta_3 ( t ) } ^2 (  q^{ \Theta ( t ) } ) ^3  \rbr*{ \tfrac{1}{3} - \tfrac{q^{\Theta ( t ) } }{4} }
	+ L \tfrac{ \abs{\Theta_1 ( t ) \Theta_3 ( t) } }{3}  (  q^{\Theta ( t ) } ) ^3   \\
	& \le (   q^{\Theta ( t ) } ) ^3 \rbr*{ - \tfrac{\abs{\Theta_3 ( t ) } ^2 }{ 12 } + \tfrac{L \abs{\Theta_3 ( t ) } }{3 } }
	= (  q^{\Theta ( t ) } ) ^3 \tfrac{\abs{\Theta_3 ( t ) } } { 12 } ( 4 L - \abs{\Theta_3 ( t ) } ) .
	\end{split}
	\end{equation}
	Hence, we obtain for all $t \in [0 , \infty )$ with $I ^{ \Theta ( t ) } \subseteq [ 0 , \frac{1}{2}  ]$ and $\abs{\Theta_3 ( t ) } \ge 4 L$ that
	$\frac{\d }{\d t } \abs{\Theta_3 ( t ) } ^2  \le 0$.
\end{cproof}

\begin{lemma}
	\label{lem:difficult:case1b}
	Consider the notation in \cref{subsec:1n:notation}
	and assume $\overline{f} < f ( 0 )$.
	Then there exists $\varepsilon  \in (0 , \nicefrac{1}{2} )$
	which satisfies for all $t \in [0, \infty )$ with $ \emptyset \not= I^{\Theta ( t ) } \subseteq [ 0 , \varepsilon ]$ and $\Theta_3 (t ) \le 0$
	that
	\begin{equation}
	\frac{\d }{\d t } \abs{\Theta_3 ( t ) } ^2  \le 0.
	\end{equation}
\end{lemma}

\begin{cproof}{lem:difficult:case1b}
	\Nobs that the fact that $f$ is continuous assures that there exists $\varepsilon  \in (0 , \nicefrac{1}{2} )$
	which satisfies for all $s \in [ 0 , \varepsilon  ]$ that $f(s ) > \overline{f}$.
	This implies for all $t \in [0, \infty )$ with $ \emptyset \not= I^{\Theta ( t ) } \subseteq [ 0 , \varepsilon  ]$
	that $\int_0^{q^{\Theta ( t ) } }  ( \overline{f} - f ( s ) ) \max \cu{\Theta_1 ( t ) s + \Theta_2 ( t ) , 0 } \, \d s \le 0$.
	Combining this with \cref{lem:integral:formulae2} demonstrates for all $t \in [0, \infty )$ with $ \emptyset \not= I^{\Theta ( t ) } \subseteq [ 0 , \varepsilon  ]$ and $\Theta_3 (t ) \le 0$
	that
	$\frac{\d }{\d t } \abs{\Theta_3 ( t ) } ^2 
	= - 2 \Theta_3 ( t ) \fG _3 ( \Theta ( t ) ) \le 0$.
\end{cproof}

\begin{lemma}
	\label{lem:difficult:case2b}
	Consider the notation in \cref{subsec:1n:notation}
	and assume $\overline{f} < f ( 0 )$.
	Then there exists $\varepsilon  \in (0 , \nicefrac{1}{2} )$
	which satisfies for all $t \in [0, \infty )$ with $ \emptyset \not= I^{\Theta ( t ) } \subseteq [ 0 ,  \varepsilon  ]$ and $\Theta_3 ( t ) > 0$
	that
	\begin{equation}
	\frac{\d }{\d t } \rbr*{ \abs{\Theta_3 ( t ) } ^2 + \tfrac{5}{8} \abs{ \Theta_1 ( t ) } ^2 } \le 0.
	\end{equation}
\end{lemma}

\begin{cproof}{lem:difficult:case2b}
	First, the fact that $f$ is continuous ensures that there exist $\beta \in (0 , \infty )$,
	$\varepsilon  \in (0 , \nicefrac{1}{2} )$
	which satisfy for all $s \in [ 0 , \varepsilon ]$
	that
	$\overline{f} + \beta < f ( s ) < \overline{f} + \frac{10 \beta }{9 }$.
	This implies for all $\theta \in \cM $ with $\emptyset \not= I^\theta \subseteq [ 0 , \varepsilon  ]$
	that
	\begin{equation}
	\int_{I^\theta } ( f ( s ) - \overline{f} ) \max \cu{\theta_1 s + \theta_2 , 0 } \, \d s \le \frac{10 \beta}{9} \int_0^{q^\theta}   \max \cu{\theta_1 s + \theta_2 , 0 } \, \d s = - \frac{5 \beta}{9} \theta_1 (   q^\theta ) ^2 
	\end{equation}
	and
	\begin{equation}
	\begin{split}
	\int_{I^\theta} ( f ( s ) - \overline{f} ) ( \theta_2^2 s - \theta_1 \theta_2 ) \, \d s 
	\ge \frac{\beta}{2} \theta_1^2 ( q^\theta ) ^2 ( 2 + (q^{\theta} ) ^2 ) .
	\end{split}
	\end{equation}
	Combining this with \cref{lem:integral:formulae2}
	demonstrates for all
	$\theta \in \cM $ with $ \emptyset \not= I^\theta \subseteq [ 0 , \varepsilon ]$
	and
	$\theta_3 \ge 0$
	that
	\begin{equation}
	\fG_3 ( \theta ) \ge \theta_1 (  q^\theta ) ^2 \rbr*{ 2 \theta_1 \theta_3   q^\theta \rbr*{ \tfrac{1}{3} - \tfrac{q^\theta}{4} } + \tfrac{10 \beta}{9} }
	\end{equation}
	and
	\begin{equation}
	\begin{split}
	\fG_1 ( \theta ) 
	&\le 2 \theta_3 \rbr*{ - \theta_1^3  \theta_3 (  q^\theta ) ^3  \rbr*{ \tfrac{1}{2} + q^\theta  } - \tfrac{\beta}{2} \theta_1^2 ( q^\theta ) ^2  ( 2 + (q^{\theta} ) ^2 ) } \\
	& = \theta_1^2 \theta_3 ( q^\theta ) ^2   \rbr*{ - \theta_1 \theta_3   q^\theta ( 1 + 2 q^\theta ) - \beta   ( 2  + ( q^\theta ) ^2 ) }.
	\end{split}
	\end{equation}
	This and the chain rule show for all 
	$t \in [0 , \infty )$
	with $ \emptyset \not= I^{\Theta ( t ) } \subseteq  [ 0 , \varepsilon ]$ and $\Theta_3 ( t ) > 0$
	that
	\begin{equation}
	\begin{split}
	&  \frac{\d }{\d t } \rbr*{ \abs{\Theta_3 ( t ) } ^2 + \tfrac{5}{8} \abs{ \Theta_1 ( t ) } ^2 } 
	= - 2 \Theta_3 ( t ) \fG_3 ( \Theta ( t ) )
	+ \tfrac{5}{4} \Theta_1 ( t ) \fG_1 ( \Theta ( t ) ) \\
	& \le 
	\Theta_1 ( t ) \Theta_3 ( t ) ( q^{ \Theta ( t ) } ) ^2 
	\Big( - \tfrac{20 \beta }{9} 
	- 4  \Theta_1 ( t ) \Theta_3 ( t ) q^{ \Theta ( t ) } \rbr*{ \tfrac{1}{3} - \tfrac{q^{ \Theta ( t ) } }{4} } \\
	& \qquad + \tfrac{5 \beta}{4} \abs{ \Theta_1 ( t ) }^2 ( 2  + ( q ^{ \Theta ( t ) } ) ^2 ) + \tfrac{5}{4} \rbr{ \Theta_1 ( t ) } ^3 \Theta_3 ( t )  q^{ \Theta ( t ) } ( 1 + 2 q^{ \Theta ( t ) } )  \Big) .
	\end{split}
	\end{equation}
	Next \nobs that 
	there exists $\eta \in (0 , \varepsilon )$
	which satisfies for all
	$t \in [0 , \infty )$
	with $ \emptyset \not= I^{\Theta ( t ) } \subseteq [0, \eta   ]$
	that
	\begin{equation}
	\tfrac{5}{4} \abs{ \Theta_1 ( t ) } ^2 ( 2  + ( q ^{ \Theta ( t ) } ) ^2 ) > \tfrac{20}{9} 
	\qqandqq - \tfrac{4}{3} + q^{ \Theta ( t ) } + \tfrac{5}{4} \rbr{ \Theta_1 ( t ) } ^2 ( 1 + 2 q^{ \Theta ( t ) } ) < 0 .
	\end{equation}
	Therefore,
	we obtain for all $t \in [0 , \infty )$
	with $ \emptyset \not= I^{\Theta ( t ) } \subseteq [ 0 , \eta  ]$ and $\Theta_3 ( t ) > 0$
	that
	\begin{equation}
	\begin{split}
	&  \frac{\d }{\d t } \rbr*{ \abs{\Theta_3 ( t ) } ^2 + \tfrac{5}{8} \abs{ \Theta_1 ( t ) } ^2 } \\
	& \le - \Theta_1 ( t ) \Theta_3 ( t ) ( q^{ \Theta ( t ) } ) ^2 \beta \rbr*{ \tfrac{20}{9} - \tfrac{5}{4} \rbr{ \Theta_1 ( t ) } ^2  ( 2  + ( q ^{ \Theta ( t ) } ) }
	\\
	& \quad + \abs{ \Theta_1 ( t ) \Theta_3 ( t ) } ^2 ( q ^{ \Theta ( t ) } ) ^3
	\rbr*{ - \tfrac{4}{3} + q^{ \Theta ( t ) } + \tfrac{5}{4} \rbr{ \Theta_1 ( t ) } ^2  ( 1 + 2 q^{ \Theta ( t ) } ) }
	\le  0.
	\end{split}
	\end{equation}
\end{cproof}

\subsection{Proof of the main boundedness result}
We now combine the results for the different cases to establish the conjecture that the entire trajectory remains bounded; see \cref{theo:bounded:1neuron} below.
The main difficulty in the proof is that the gradient flow may change between the different regimes.

\begin{theorem}
	\label{theo:bounded:1neuron}
	Consider the notation in \cref{subsec:1n:notation} and assume that $f$ is Lipschitz continuous. Then $\sup_{t \in [0 , \infty ) } \norm{\Theta(t) } < \infty $.
\end{theorem}

\begin{remark}
	The assumption that $f$ is Lipschitz is only needed in the special cases $\overline{f} = f(0)$ (see \cref{lem:endpoint:mean:value2})
	and
	$\overline{f} = f(1)$
	(see \cref{lem:endpoint:mean:value}).
	If one assumes $f(0) \not= \overline{f} \not= f(1)$ it is sufficient if $f$ is merely continuous.
\end{remark}

\begin{cproof}{theo:bounded:1neuron}
	First \nobs that if there exists $t \in [0 , \infty )$ with $\mu ( I^{ \Theta ( t ) }  ) = 0$ then $\fG ( \Theta ( t ) ) = 0$. By uniqueness of solutions (since $\fG$ is locally Lipschitz on $\cM$), we obtain for all $u \in [0 , \infty )$ that $\fG ( \Theta ( u ) ) = 0$ and, hence, $\Theta ( u ) = \Theta ( 0 )$. In this case the statement clearly holds.
	
	From now on we assume $\forall \, t \in [0 , \infty ) \colon \mu ( I^{ \Theta ( t ) } ) > 0 $.
	We consider the case $\overline{f} < \min \cu{f(0) , f(1) }$. The remaining cases are analogous, using \cref{lem:endpoint:mean:value,lem:endpoint:mean:value2}.
	\Nobs that \cref{lem:difficult:case1,lem:difficult:case2,lem:difficult:case1b,lem:difficult:case2b} assure that there exists $\varepsilon \in (0 , \nicefrac{1}{2} )$ which satisfies the following properties:
	\begin{enumerate} [label = (\Roman*)]
		\item
		\label{theo:bounded:proof:item1}
		It holds for all $t \in [0, \infty )$ with $ \emptyset \not= I^{\Theta ( t ) } \subseteq [ 1 - \varepsilon , 1 ]$ and $\Theta_3 (t ) \le 0$
		that $\frac{\d}{\d t} ( \abs{ \Theta_3 ( t ) } ^2 ) \le 0$,
		\item
		\label{theo:bounded:proof:item2}
		it holds for all $t \in [0, \infty )$ with $ \emptyset \not= I^{\Theta ( t ) } \subseteq [ 0 , \varepsilon ]$ and $\Theta_3 (t ) \le 0$
		that $\frac{\d}{\d t} ( \abs{ \Theta_3 ( t ) } ^2 ) \le 0$,
		\item
		\label{theo:bounded:proof:item3}
		it holds for all $t \in [0, \infty )$ with $ \emptyset \not= I^{\Theta ( t ) } \subseteq [ 1 - \varepsilon , 1 ]$ and $\Theta_3 ( t ) > 0$
		that
		$\frac{\d }{\d t } \rbr*{ \abs{\Theta_3 ( t ) } ^2 - \tfrac{5}{8 } \abs{ \Theta_1 ( t ) - 2^{- \nicefrac{1}{2} } } ^2 } \le 0$,
		and
		\item
		\label{theo:bounded:proof:item4}
		it holds for all for all $t \in [0, \infty )$ with $ \emptyset \not= I^{\Theta ( t ) } \subseteq [ 0 ,  \varepsilon  ]$ and $\Theta_3 ( t ) > 0$
		that $\frac{\d }{\d t } \rbr*{ \abs{\Theta_3 ( t ) } ^2 + \tfrac{5}{8} \abs{ \Theta_1 ( t ) }^2 } \le 0$.
	\end{enumerate}
	Next let $\cT \subseteq [0 , \infty )$
	satisfy
	\begin{equation}
	\cT = \cu{t \in [0 , \infty ) \colon \mu ( I ^{ \Theta ( t ) } ) \in [ \varepsilon , 1 )  }
	\cup \cu{t \in [0 , \infty ) \colon \abs{\Theta_1 ( t ) } \ge \varepsilon , \, \mu ( I^{\Theta ( t ) } ) = 1 }.
	\end{equation}
	\Nobs that \cref{lem:bounded:case:simple1,lem:bounded:case:simple2} imply that $\fC = \abs{\Theta_3 ( 0 ) } ^2 + \sup_{t \in \cT } \abs{\Theta_3 ( t ) }^2  < \infty$.
	Now let $\tau \in [0 , \infty )$ be arbitrary, we will show that $\abs{\Theta_3 ( \tau ) }^2  < \fC + 3$.
	Define
	\begin{equation}
	\label{theo:bounded:proof:eq:defu}
	u = \sup \rbr*{  \cu{t \in [0 , \tau) \colon \abs{\Theta_3 ( t ) } ^2 \le \fC  } }
	\end{equation}
	and assume without loss of generality that $u < \tau$
	and $\tau \notin \cT$.
	\Nobs that this implies that $\mu ( I^{ \Theta ( \tau )}) < \varepsilon$ or $\mu ( I ^{ \Theta ( \tau ) } ) = 1$.
	We now consider four cases.
	
	\setcounter{case}{0}
	\begin{case}
		Assume $\mu ( I^{ \Theta ( \tau ) } ) = 1$. In this case, we necessarily have $\forall \, t \in (u , \tau ) \colon \mu ( I^{ \Theta ( t ) } ) = 1$.
		Indeed, otherwise by continuity of $t \mapsto \mu ( I^{ \Theta ( t ) } )$ there would exist $t \in (u , \tau )$ with $\mu ( I^{ \Theta ( t ) } ) \in ( \varepsilon , 1 )$. Hence $t \in \cT$ and $\abs{\Theta_3 ( t ) } ^2 \le \fC$, which contradicts \cref{theo:bounded:proof:eq:defu}.
		Furthermore, from \cref{theo:bounded:proof:eq:defu} we obtain for all $t \in (u , \tau )$ that $\abs{\Theta_1 ( t ) } < \varepsilon $. In addition, \cref{prop:neuron:active:bounded} ensures for all $t \in ( u , \tau )$
		that $\frac{\d}{\d t } \rbr*{ \abs{\Theta_3 ( t ) } ^2 + \ln ( 1 - \abs { \Theta_1 ( t ) } ^2 ) } = 0$.
		Hence, we obtain that
		\begin{equation}
		\begin{split}
		\abs{\Theta_3 ( \tau ) } ^2 
		&\le \abs{\Theta_3 ( u ) } ^2 + \abs{ \ln ( 1 - \abs { \Theta_1 ( u ) } ^2 ) - \ln ( 1 - \abs { \Theta_1 ( \tau ) } ^2 ) } \\
		&\le \abs{\Theta_3 ( u ) } ^2 + 2 \abs{ \ln ( 1 - \varepsilon ^2 ) } \le \fC  + 2 \abs{ \ln ( \tfrac{3}{4} ) } < \fC + 3.
		\end{split}
		\end{equation}
	\end{case}
	\begin{case}
		Assume $\mu( I^{ \Theta ( \tau ) } ) < \varepsilon$ and $\Theta_3 ( \tau ) < 0$.
		Since $[0 , \infty ) \ni t \mapsto \Theta_3 ( t ) \in \R$ and
		$ [0 , \infty ) \ni  t \mapsto \mu ( I^{ \Theta ( t ) } ) \in \R$ are continuous,
		\cref{theo:bounded:proof:eq:defu} shows for all $t \in ( u , \tau )$ that $\Theta_3 ( t ) < - \sqrt{\fC} < 0$ and $( I ^{ \Theta ( t ) } \subseteq [ 0 , \varepsilon ] ) \vee ( I ^{ \Theta ( t ) } \subseteq [ 1 - \varepsilon , 1 ] )$. 
		\ref{theo:bounded:proof:item1} and \ref{theo:bounded:proof:item2}
		therefore imply for all $t \in ( u , \tau )$ that $\frac{\d}{\d t} ( \abs{ \Theta_3 ( t ) } ^2 ) \le 0$.
		Hence, we obtain that $\abs{\Theta_3 ( \tau ) } ^2 \le \abs{\Theta_3 ( u ) } ^2 \le \fC$.
	\end{case}

	\begin{case}
		Assume $\mu( I^{ \Theta ( \tau ) } ) < \varepsilon$, $\Theta_3 ( \tau ) > 0$,
		and $I^{ \Theta ( \tau ) } \subseteq [1 - \varepsilon , 1 ]$.
		By continuity of $ [0 , \infty ) \ni  t \mapsto \Theta_3 ( t ) \in \R$ and \cref{theo:bounded:proof:eq:defu} we obtain for all $t \in ( u , \tau )$ that $\Theta_3 ( t ) > \sqrt{\fC} > 0$ and $I ^{ \Theta ( t ) } \subseteq [1 - \varepsilon , 1 ] $. \ref{theo:bounded:proof:item3} therefore demonstrates for all $t \in (u , \tau )$ that $\frac{\d }{\d t } \rbr*{ \abs{\Theta_3 ( t ) } ^2 - \tfrac{5}{8} \abs{ \Theta_1 ( t ) - 2^{ - \nicefrac{1}{2} } } ^2 } \le 0$.
		This yields that
		\begin{equation}
		\begin{split}
		\abs{\Theta_3 ( \tau ) } ^2 
		&\le \abs{\Theta_3 ( u ) } ^2 + \tfrac{5}{8} \abs{  \abs{ \Theta_1 ( u ) - 2^{ - \nicefrac{1}{2} } } ^2  - \abs{ \Theta_1 ( \tau ) - 2^{ - \nicefrac{1}{2} } } ^2 } \\
		&\le \fC  + \tfrac{5}{2} < \fC + 3 .
		\end{split}
		\end{equation}
	\end{case}
	
	\begin{case}
		Assume $\mu( I^{ \Theta ( \tau ) } ) < \varepsilon$, $\Theta_3 ( \tau ) > 0$,
		and $I^{ \Theta ( \tau ) } \subseteq [0 , \varepsilon ]$. 
		By continuity of $ [0 , \infty ) \ni  t \mapsto \Theta_3 ( t ) \in \R $ and \cref{theo:bounded:proof:eq:defu} we obtain for all $t \in ( u , \tau )$ that $\Theta_3 ( t ) > \sqrt{\fC} > 0$ and $I ^{ \Theta ( t ) } \subseteq [0 , \varepsilon ] $. \ref{theo:bounded:proof:item4} therefore proves for all $t \in (u , \tau )$ that $\frac{\d }{\d t } \rbr*{ \abs{\Theta_3 ( t ) } ^2 + \tfrac{11}{18} \abs{ \Theta_1 ( t ) }^2 } \le 0$.
		This implies that
		\begin{equation}
		\begin{split}
		\abs{\Theta_3 ( \tau ) } ^2 
		&\le \abs{\Theta_3 ( u ) } ^2 + \tfrac{5}{8} \abs*{ \abs{ \Theta_1 ( u ) }^2 - \abs{ \Theta_1 ( \tau ) } ^2 } \le \fC  + \tfrac{5}{4} < \fC + 3 .
		\end{split}
		\end{equation}
	\end{case}
\end{cproof}

\subsection*{Acknowledgments}
The second and third authors acknowledge funding by the Deutsche Forschungsgemeinschaft 
(DFG, German Research Foundation) under Germany's Excellence Strategy 
EXC 2044-390685587, Mathematics M\"{u}nster: Dynamics-Geometry-Structure.
This project has been partially supported by the startup fund project of Shenzhen Research Institute of Big Data under grant No. T00120220001.

\end{document}